\documentclass[a4paper,10pt]{article}
\usepackage{setspace}
\usepackage[utf8]{inputenc}
\usepackage{amsmath}
\usepackage{amssymb}
\usepackage{amsfonts}
\usepackage{stmaryrd} 
\usepackage{dsfont} 

\usepackage{import}
\usepackage{todonotes}

\usepackage{tabularx, multirow} 
\usepackage{threeparttable}
\usepackage[]{booktabs}         

\usepackage{float}
\usepackage{siunitx}

\usepackage{subfigure}

\usepackage{enumerate}

\usepackage[mathlines]{lineno}
\newcommand*\patchAmsMathEnvironmentForLineno[1]{%
  \expandafter\let\csname old#1\expandafter\endcsname\csname #1\endcsname
  \expandafter\let\csname oldend#1\expandafter\endcsname\csname end#1\endcsname
  \renewenvironment{#1}%
     {\linenomath\csname old#1\endcsname}%
     {\csname oldend#1\endcsname\endlinenomath}}%
\newcommand*\patchBothAmsMathEnvironmentsForLineno[1]{%
  \patchAmsMathEnvironmentForLineno{#1}%
  \patchAmsMathEnvironmentForLineno{#1*}}%
\AtBeginDocument{%
\patchBothAmsMathEnvironmentsForLineno{equation}%
\patchBothAmsMathEnvironmentsForLineno{align}%
\patchBothAmsMathEnvironmentsForLineno{flalign}%
\patchBothAmsMathEnvironmentsForLineno{alignat}%
\patchBothAmsMathEnvironmentsForLineno{gather}%
\patchBothAmsMathEnvironmentsForLineno{multline}%
}

\newcommand{\lcap}{\,l_{\text{cap}}}

\newcommand{\hjur}{h_{\text{Jurin,2D}}}

\newcommand{\fcL}{\lambda_{\Gamma}}
\newcommand{\cl}{\Gamma}

\usepackage{authblk}
\usepackage{microtype}
\author[1]{Gründing, D.}
\author[2]{Smuda, M.}
\author[3]{Antritter, T.}
\author[1]{Fricke, M.}
\author[4]{Rettenmaier, D.}
\author[2]{Kummer, F.}
\author[3]{Stephan, P.}
\author[1]{Marschall, H.}
\author[1]{Bothe, D.}
\affil[1]{Institute for Mathematical Modeling and Analysis, Technische Universit\"at Darmstadt}
\affil[2]{Institute for Fluid Dynamics, Technische Universit\"at Darmstadt}
\affil[3]{Institute for Technical Thermodynamics, Technische Universit\"at Darmstadt}
\affil[4]{Institute for Fluid Mechanics and Aerodynamics, Technische Universit\"at Darmstadt}

\title{Capillary Rise - A Computational Benchmark for Wetting Processes}

\doublespace
\begin{document} 
 
\maketitle
 

\begin{abstract}
Four different numerical approaches are compared for the rise of liquid between two parallel plates. These are an  Arbitrary Lagrangian-Eulerian method (OpenFOAM solver interTrackFoam), a geometric volume of fluid code (FS3D), an algebraic volume of fluid method (OpenFOAM solver interFoam), and a level set approach (BoSSS). The first three approaches discretize the bulk equation using a finite volume method while the last one employs an extended discontinuous Galerkin discretization. 
The results are compared to ODE models which are the classical rise model and an extended model that incorporates a Navier slip boundary condition on the capillary walls and levels at a corrected stationary rise height.
All physical parameters are based on common requirements for the initial conditions, short simulation time, and a non-dimensional parameter study. The comparison shows excellent agreement between the different implementations with minor quantitative deviations for the adapted interFoam implementation.  
While the qualitative agreement between the full solutions of the continuum mechanical approach and the reference model is good, the quantitative comparison is only reasonable, especially for cases with increasing oscillations. Furthermore, reducing the slip length changes the solution qualitatively as oscillations are completely damped in contrast to the solution of the ODE models. 
To provide reference data for a full continuum simulation of the capillary rise problem, all results are made available online.
\end{abstract} 

\textbf{Keywords:} wetting, capillary rise, benchmark comparison, direct numerical simulation, Navier slip, multiphase flow

\newpage
\textbf{Highlights}
\begin{itemize}
 \item non-dimensional comparison between a simplified reference model and solutions of the full PDE system by four different numerical approaches/codes for numerical slip and Navier slip boundary conditions
 \item using ``numerical slip`` dampens rise dynamics with increasing mesh resolution and renders results unusable while applying a Navier slip boundary condition shows mesh convergent results
 \item reducing the slip length in the full numerical solution changes the behavior qualitatively in comparison to the reference solution as rise height oscillations disappear
 \item benchmark data available online
\end{itemize}

\newpage

\section{Introduction}\label{sec:introduction}
The problem of a liquid penetrating into spaces smaller than the capillary length can be found in various processes. The combination of multiple capillaries forms a basic model for liquid movement in porous media described, e.g., by network models. Hence the penetration of liquid into cavities is one of the basic physical phenomena relevant for various industrial applications such as crude oil recovery from oil sand, or the penetration of dye into paper in printing applications.

Due to its geometric simplicity and experimental accessibility, the rise of liquid in a single capillary has also been subject to investigation by a large number of researchers over the last century. While the stationary rise height of a liquid in a cylindrical capillary is available since \cite{Jurin1719}, one of the first quantitative models of a liquid rising in a cylindrical capillary over time was given by \cite{Lucas1918,Washburn1921}. The model consists of an ad-hoc force balance accounting for inertia, viscous, gravitation and surface tension forces. It was used to model not only the rise of liquid in a single capillary, but also the penetration of liquid into porous media such as filtration paper. This basic model has been extended in various ways. In \cite{Szekely1971}, the influence of the reservoir feeding liquid into the capillary is considered. Their adaptation allows to remove the singularity of the classical model that occurs for vanishing initial height. Improving upon the work of \cite{Szekely1971}, the influence of the reservoir is analyzed in detail in \cite{Levine1976} yielding a similar term for regularization, and a compelling argument for the pressure influence at the inflow of the capillary. 
 
A dimensional analysis of the classical model is provided in \cite{Fries2009}. There, a single non-dimensional group is identified. Furthermore, it is demonstrated how various analytical solutions from literature correspond to limits of an identified non-dimensional group. Estimates for the different rise regimes corresponding to the limiting cases are provided by \cite{Quere1997}. This includes an approximate condition for the occurrence of rise height oscillations. In \cite{Stange2003}, the classical model is extended to incorporate various effects of a cylindrical reservoir. These include the free surface of the reservoir, and a model for the pressure field in the inflow region of the capillary. In addition, the decay of interface oscillations was added to the model.

Using the approximate analytic description of \cite{Huh1971} for a velocity field in the wedge near the contact line, \cite{Levine1980} provides local velocity fields for the quasi-stationary rise of a liquid in a capillary. The author also gives an improved prediction of the stationary height.

The problem of a liquid rising in a capillary is interesting as a benchmark for various reasons. First, the general setup, can be simplified to a rise of a liquid between two parallel plates infinitely extended along the median plane. This reduces the necessary run time due to less degrees of freedom in comparison to the full 3D geometry. Furthermore, various different models for the rise of a liquid in a capillary are available. This model family does also contain extensions to incorporate various geometric and physical influences. Finally, the rise of a liquid in a capillary is comparably well accessible in physical experiments. While for stationary solutions of wetting problems non-trivial interface geometries can be obtained , analytic or approximate reference solutions for a moving contact line are rather scarce. Hence, we choose the rise of a liquid between to parallel plates as a benchmark case to assess the capability of four different CFD codes to simulate this problem.

Surprisingly few authors deal with a solution of the full continuum mechanical description of the rise of a liquid in a capillary. In \cite{Sprittles2012} the quasi-stationary propagation of a liquid-gas meniscus through a cylindrical tube is investigated where the slip length of $\SI{1}{\nano \meter}$ is resolved. Results are compared to local asymptotic solutions of the velocity field near the contact line. The case of forced liquid flow into a cylindrical capillary is analyzed in \cite{Sui2011} using a rotationally symmetric level set approach. The resolved flow is governed by a prescribed Poiseuille flow velocity field at the inflow of the capillary. Among others, a quasi-stationary, a wavy, and a jetting regime are identified. In \cite{Gruending2019}, a model for liquid rising between two planar plates is derived for the case, where a slip boundary condition is applied on the complete surface of the capillary. It is demonstrated, that for slip on the scale of the capillary, mesh-convergent results can be obtained that show a qualitatively good agreement with predictions of the adapted rise model. Furthermore, a correction for the stationary rise height is introduced to the considered class of ODE models.
 
There is no publication known to the authors concerning a detailed analysis for the full continuum mechanical description of liquid rising in a capillary up to its stationary height for a larger range of parameters. The present work intends to provide benchmark results for a surface tension driven test case solved with four different numerical approaches and a comparison with available models from the literature. Due to reduced computational costs, 2D or rotational symmetric geometries are preferable for the purpose of a benchmark case.

\section{Mathematical model for a liquid rising in a gap}
\begin{figure}[H]
 \begin{center}
 \def\svgwidth{.3\textwidth} 
 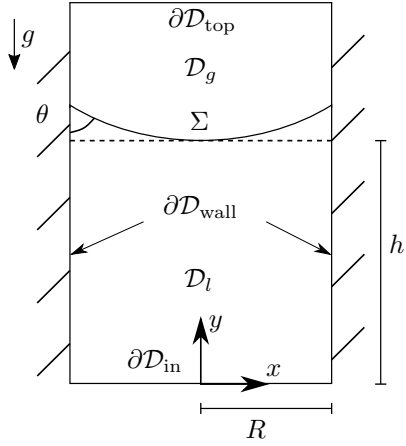 
 \caption{Domain and boundary names for liquid rising in a capillary.} 
  \label{fig:capRiseDomain} 
 \end{center} 
\end{figure}

Various models for the rise of liquids in capillaries are available. These typically give a prediction of the apex height $h$ over time; see Figure~\ref{fig:capRiseDomain}, where the complete domain $\mathcal{D}$ consists of a liquid part ($l$) and a gaseous part ($g$) separated by the interface $\Sigma(t)$. The domain is bounded by $\partial \mathcal{D}_{\text{in}}$ at the bottom and the outflow boundary $\partial \mathcal{D}_{\text{top}}$ at the upper part of the capillary. Both boundaries need to allow fluid flowing through them. The capillary walls are denoted by $\partial \mathcal{D}_{\text{wall}}$. 

For small E\"otv\"os numbers $\text{Eo} = \Delta \rho g R^2 / \sigma$, surface tension forces dominate gravitational forces and the interface can be well approximated by a circular section \cite{Concus1968}. In addition, for typical gas-liquid systems, the influence of the gas phase can be neglected. Under these assumptions, the stationary rise height $h_{\infty}^{\text{apex}}$ for a liquid in a capillary gap can be calculated via
\begin{align}
 h_{\infty}^{\text{apex}}
 &=
 \hjur - \hat{h},\qquad 
 \hjur 
 =
 \frac{\sigma \cos \theta}{R\rho g}, \qquad 
 \hat{h}
 =
\frac{R}{2\cos\theta} \left(2 - \sin \theta - \frac{\arcsin \cos \theta}{\cos \theta}\right),
 \label{eq:stationary_rise_height}
\end{align}
where $h_{\text{Jurin}}$ is the estimated rise height of liquid in a gap by the popular formula of \cite{Jurin1719}. The height correction $\hat{h}$ corresponds to incorporating the liquid volume in the interface region, i.e.  the area between the dashed apex height and solid interface line in Figure~\ref{fig:capRiseDomain}, cf. \cite{Gruending2019}.

\subsection{Continuum mechanical model}\label{subsec:govEqs}
\newcommand{\jump}[1]{\left\llbracket #1 \right\rrbracket}
\newcommand{\nsigma}{\boldsymbol{n}_\Sigma}
\newcommand{\domain}{\mathcal{D}}
\newcommand{\nouter}{\boldsymbol{n}_{\partial\mathcal{D}}}

We briefly recall the sharp interface two-phase Navier Stokes model (see, e.g., \cite{Pruss.2016},  \cite{Slattery.1999}). Assuming incompressible Newtonian fluids, the local balance equations for mass and momentum in the bulk phases read
\begin{align}
\nabla \cdot \boldsymbol{u} = 0 \quad \text{in} \quad \domain\setminus\Sigma(t),\label{eq:massCons}\\
\rho (\partial_t \boldsymbol{u} + \boldsymbol{u} \cdot \nabla \boldsymbol{u}) - \mu \Delta \boldsymbol{u} + \nabla p = \boldsymbol{f} \quad \text{in} \quad \domain\setminus\Sigma(t), \label{eq:momCons}
\end{align}
where $\mu$ denotes the dynamic viscosity and $\boldsymbol{f}$ is an external body force density. In this paper we only consider gravity leading to $\boldsymbol{f} = \rho\boldsymbol{g}$. If no mass is transferred across the fluid-fluid interface $\Sigma(t)$, the normal component of the bulk velocities are continuous and equal to the normal velocity $V_\Sigma$ of $\Sigma$, i.e.
\begin{align}
\boldsymbol{u}^+ \cdot \nsigma = \boldsymbol{u}^- \cdot \nsigma \quad \Leftrightarrow \quad \jump{\boldsymbol{u}} \cdot \nsigma = 0 \quad &\text{on} \quad \Sigma(t),\\
V_\Sigma = \boldsymbol{u}^\pm \cdot \nsigma \quad &\text{on} \quad \Sigma(t).\label{eq:kinematic_condition}
\end{align}
Moreover, it is assumed that there is no slip between the fluid phases. In this case the tangential velocity components are also continuous leading to
\begin{align}
\label{eq:jump_condition_velocity}
\jump{\boldsymbol{u}} = 0 \quad \text{on} \quad \Sigma(t).
\end{align}
Above, the symbol $\jump{\cdot}$ denotes the jump of a quantity across the fluid-fluid interface, i.e.
\[ \jump{\phi}(t,\boldsymbol{x}) := \lim_{h\rightarrow 0^+} \phi(t,\boldsymbol{x}+h\nsigma)-\phi(t,\boldsymbol{x}-h\nsigma). \]
The effect of interfacial tension is mediated via the transmission condition for the stress. Assuming constant surface tension it reads
\begin{align}\label{eq:jump_condition_momentum}
\jump{p \mathds{1} - \mu(\nabla \boldsymbol{u} + \nabla \boldsymbol{u}^{\sf T})} \nsigma = \sigma \kappa \nsigma  \quad \text{on} \quad \Sigma(t),
\end{align}
where $p$ denotes the pressure, $\kappa := - \nabla_\Sigma \cdot \nsigma$ is the mean curvature of the interface and $\sigma > 0$ denotes the surface tension coefficient. Below, we use the symbol for the viscous stress tensor $\boldsymbol{S}$, i.e.
\[ \boldsymbol{S}= \mu(\nabla \boldsymbol{u} + \nabla \boldsymbol{u}^{\sf T}). \]

\paragraph{Boundary conditions at the solid boundary:}It is a well-known fact that the no slip condition (expressed in a frame of reference where the solid wall is at rest), i.e.
\begin{align}
\boldsymbol{u} = 0,
\end{align}
is incompatible with a moving contact line \cite{Huh1971}. Hence, some mechanism is required allowing for a non-zero tangential component of the interface velocity. In the present paper we consider the Navier slip boundary condition. It can be derived by thermodynamic considerations and takes the form
\begin{align}
\label{eq:navier_condition}
\boldsymbol{u} \cdot \nouter = 0, \quad \lambda \boldsymbol{u}_\parallel + (\boldsymbol{S}\nouter)_\parallel = 0 \quad \text{on} \quad \partial\domain_{\text{wall}},
\end{align}
where $\nouter$ denotes the outer unit normal field to $\domain$. The model parameter $\lambda > 0$ is a friction coefficient determining the amount of slip. By dividing equation \eqref{eq:navier_condition} by $\lambda$, we obtain the equivalent formulation
\begin{align}
\label{eq:navier_condition_v2}
\boldsymbol{u} \cdot \nouter = 0, \quad \boldsymbol{u}_\parallel + L \, \frac{\partial \boldsymbol{u}_\parallel}{\partial \nouter} = 0 \quad \text{on} \quad \partial\domain_{\text{wall}}.
\end{align}
The quantity $L=\mu/\lambda$ has the dimension of a length and is called \emph{slip length}. Note that for $L=0$ we recover the no slip condition, while the limit $L \rightarrow \infty$ is also considered in the literature under the name ``free slip'' or ``perfect slip'' condition.
\paragraph{Boundary condition at the contact line:} To complete the model, one has to specify a boundary condition for the contact angle. In the present study we consider the case of a fixed contact angle equal to the equilibrium value, i.e.
\begin{align}
\label{eq:contact_angle_boundary_condition}
\theta = \theta_{\text{e}} \quad \text{at} \quad \Gamma.
\end{align}

\paragraph{Boundary conditions at the inflow and outflow boundaries:}
In addition to the boundary conditions at the physical boundaries, boundary conditions at the artificial inflow and outflow boundaries have to be specified. In the present study we consider the boundary conditions (see \cite{Gruending2019b} for a discussion of the pressure boundary condition)

\begin{align}
 \partial_n u_n = 0, \quad p = 0 \quad &\text{on} ~ \partial \domain_{\text{in}} \cup \partial \domain_{\text{top}}
\end{align}
together with
\[ \rho = \rho_l \quad \text{on} \quad \domain_{\text{in}}, \quad \rho = \rho_g \quad \text{on} \quad \domain_{\text{top}}. \]

\subsection{Reference models for the capillary rise problem}\label{sub:ode_model}
Various simplified models for the capillary rise problem are available in the literature. The classical description goes back to \cite{Lucas1918,Washburn1921}. Adapted to the flow between two planar plates it reads
\begin{align}\label{eq:ode_classical}
 \rho \frac{\text{d} }{\text{d}t}(\dot{h}h) 
 &=
 -\frac{3\mu}{R^2}\dot{h}h - \rho g h  +  \frac{\sigma \cos \theta_e}{R}, \quad h(t_0)>0,\, \dot{h}(t_0)=\dot{h}_0.
\end{align}
The model is based on a force balance of inertia, viscous, gravitation and surface tension forces. The viscous forces are based on the assumption of a Poiseuille flow with no-slip boundary conditions on the capillary wall. Furthermore, the stationary rise height of the classical model is consistent with Jurin's height $\hjur$. 

As will be demonstrated in Section \ref{sec:setups_results}, no-slip boundary conditions may yield non-convergent rise height behavior. Navier slip boundary conditions are used to circumvent this problem. Physical slip lengths for the Navier slip boundary condition are typically quantified on the scale of nanometers \cite{Qian2003}. As we do not seek a direct comparison to experiments but rather a code to code comparison, it is sufficient to use a slip length on the scale of the capillary. However, the classical model does not incorporate any influence of a slip boundary condition. This influence is incorporated by an extended model due to \cite{Gruending2019} which reads
\begin{align}\label{eq:ode_extended}
 \rho \frac{\text{d} }{\text{d}t} (\dot{h}(h  + \hat{h})) 
 =&-\frac{3\mu\dot{h}}{R(R+3L)}(h+\hat h)
- \rho g (h + \hat{h} ) \nonumber \\ 
&+ \frac{\sigma \cos \theta_e}{R}  
 +\rho \dot{h}^2\frac{3(15L^2 + 10LR + 2R^2)}{5(R+3L)^2},
 \quad h(t_0)\geq 0,\,\dot{h}(t_0)=\dot{h}_0.
\end{align}
The model also includes the influence of mass between a cross section at apex height and the interface. Firstly, this yields the regularization in the inertia term. Secondly, the extended model is consistent with the stationary rise height predicted by \eqref{eq:stationary_rise_height} in contrast to the classical model \eqref{eq:ode_classical}. It has been shown in \cite{Gruending2019} that the difference between Jurin's height and the stationary height (see \eqref{eq:stationary_rise_height}) is relevant for a comparison between a rise model such as \eqref{eq:ode_classical} and a full solution of the continuum mechanical problem \cite{Gruending2019}. Moreover, the extended model includes the influence of a Navier slip boundary condition on the capillary walls in the viscous term and the convection term ($\dot{h}^2$ term). Hence, the problem to resolve both, the capillary scale and the nanometer scale for realistic slip lengths, is significantly relieved.

The assumptions necessary to obtain the reference models from the continuum mechanical mass and momentum conservation equations are outlined in, e.g., \cite{Gruending2019}. The derivation assumes that the interface has the shape of a spherical cap for the case between two plates (spherical cap for a cylinder). This assumption requires a constant contact angle. Furthermore, only the velocity component in rise direction is considered, neglecting any influence of the viscous dissipation in contact line vicinity. At the inflow boundary, a Poiseuille flow profile is assumed yielding the $\dot{h}^2$ term in contrast to the classical model. This assumption is combined with a constant pressure at the inflow boundary. For $L=0$ and $\hat h=0$, the extended model reduces to the classical description if the convection term is neglected.

\subsubsection{Scaling of the approximate model}\label{sub:scaling}
For the classical capillary rise model without slip, a dimensional analysis is given in \cite{Fries2009} and a single relevant dimensional group $\Omega$ has been identified. The various limits for $\Omega$ correspond the special cases for which analytic solutions are available. An analog non-dimensionalization has been performed for the model \eqref{eq:ode_extended}. The relevant factors are listed in Table~\ref{tab:scalingCoeffs}.  The coefficients for the 2D case have been defined analog to the 3D case from \cite{Fries2009}. 

\begin{table}[H] 
\caption{Scaling coefficients for 2D/3D}
\label{tab:scalingCoeffs}
\centering
  \begin{tabular}{ l  l  l  l  l}
   & a & b & c & $\Omega$\\ \hline 
    2D &$\rho R/ (\sigma \cos \theta_e)$ & $3 \mu/ (R\sigma \cos \theta_e)$ & $\rho g R/(\sigma \cos \theta_e)$ & $\sqrt{9\sigma \cos \theta_e \mu^2/(\rho^3 g^2 R^5)}$\\ 
    3D &$\rho R/ (2\sigma \cos \theta_e)$ & $4 \mu/ (R\sigma \cos \theta_e)$ & $\rho g R/(2\sigma \cos \theta_e)$ & $\sqrt{128\sigma \cos \theta_e \mu^2/(\rho^3 g^2 R^5)}$\\
  \end{tabular}
\end{table}

Choosing the coefficients $a$, $b$ and $c$ as in table \ref{tab:scalingCoeffs}, the time and length scales as well as the non-dimensional group $\Omega$ can be expressed in the same way as in \cite{Fries2009}. This gives rise to three different scalings denoted by I, II, and III with scalings given in Table \ref{tab:timeAndLengthScales}. Considering, for example scaling I, the dimensionless time $t^*$ is defined by $t^*=c^2/b t$ and the dimensionless rise height $h^*$ by $h^*=ch$.
From here on we will use dimensionless quantities where applicable and drop the ``$^*$''.
\begin{table}[H]
\caption{Time and length scales}
\label{tab:timeAndLengthScales}
\centering
  \begin{tabular}{ l l l l }
   & I & II & III \\ \hline
  t & $c^2/b$ & $\sqrt{c^2/a}$ & $b/a$\\
  h & $c$ & $c$ & $b/\sqrt{2a}$ 
  \end{tabular}
\end{table}
With the quantities listed in Table~\ref{tab:timeAndLengthScales}, the non-dimensional group $\Omega$ can be expressed for both, the 2D and the 3D case, by
\begin{align}
 \Omega = \sqrt{\frac{b^2}{ac^2}}.
\end{align}
To describe the influence of the slip length in the 2D case, we define
\begin{equation}
 S :=\frac{L}{R}, \quad 
 K
 :=
 \frac{1}{1+3S}, \quad Q := \frac{3(15S^2 + 10S+2)}{5(1+3S)^2}.
\end{equation}
Using the three different scalings listed in Table \ref{tab:scalingCoeffs}, the rise model \eqref{eq:ode_extended} is non-dimensionalized. The scaled model for all three scalings is listed in Table \ref{tab:scaledModels}.
\begin{table}[H]
\caption{Scaled extended model with Navier slip b.c.}
\label{tab:scaledModels}
\centering
  \begin{tabular}{ l  l}
  Scaling & Scaled Equation \\
  \hline
    I   & $1/\Omega^2 \partial_t(\dot{h}(h+\hat{h})) + K \dot{h}(h+\hat{h}) + h+\hat{h} = 1 + Q/\Omega^2 \dot{h}^2 $ \\ 
    II  & $\partial_t(\dot{h}(h+\hat{h})) + K \Omega (h+\hat{h}) \dot{h} + h + \hat{h} = 1 + Q \dot{h}^2$ \\
    III & $2 \partial_t(\dot{h}(h+\hat{h})) + 2 K \dot{h}(h+\hat{h}) + \sqrt{2}/\Omega (h+\hat{h}) = 1 + 2Q\dot{h}^2$ 
  \end{tabular}
\end{table}  
Note that for $L=0$ it follows that $K=1$. Hence, the  scaling for the model with Navier slip is consistent with the scaling for the model with no slip boundary conditions \cite{Fries2009,Gruending2019}.

\section{Numerical methods}\label{sec:numerical_methods}
In the following we give a brief overview regarding the main aspects of each participating numerical method. These are: 1. an Arbitrary Lagrangian Eulerian approach (OpenFOAM solver interTrackFoam), 2. a geometric volume of fluid code (FS3D), 3. an algebraic volume of fluid method (OpenFOAM solver interFoam) and 4. level set implementation (BoSSS).
 
  \subsection{Arbitrary Lagrangian Eulerian method - interTrackFoam}
One possible approach to solve the analytic description of the wetting problem given in Section \ref{subsec:govEqs} is an Arbitrary Lagrangian Eulerian (ALE) method. With this approach, the simulation domain is tessellated with non-intersecting control volumes $V(t) \subset \Omega(t), \mathring{V}(t) \cap \Sigma(t)=\emptyset$ such that faces of the mesh represent the discretized interface. Integrating the mass and momentum conservation equations \eqref{eq:massCons} and \eqref{eq:momCons} over such a control volume, yields the integral equations
\newcommand{\dx}[1]{\,\mathrm{d}{#1}}
\newcommand{\vel}{\boldsymbol{u}}
\newcommand{\normal}{\boldsymbol{n}}
\newcommand{\ddt}{\frac{\text{d}}{\text{d}t}}
\begin{align}
\ddt \int_{V(t)} \rho \dx{V} &= \int_{\partial V(t)} (\vel_b - \vel) \cdot \normal \dx{o}\label{eq:integral_mass_conservation} \\
\ddt \int_{V(t)} \rho \vel \dx{V}
 &= 
 \int_{\partial V(t)} \rho \vel \otimes (\vel_b - \vel) \normal \dx{o}+ 
 \int_{\partial V(t)} \boldsymbol{S} \normal \dx{o} + \int_{V(t)} \rho \mathbf{g} \dx{V} \label{eq:integral_momentum_conservation}
\end{align}
where $\vel_b$ is the velocity of the control volume boundary $\partial V(t)$. For systems such as a liquid (l) displacing gas (g) from a capillary, $\mu_g/\mu_l\ll 1$ and $\rho_g/\rho_l\ll 1$. This allows to neglect the influence of the gas on the liquid phase reducing the transmission conditions \eqref{eq:jump_condition_momentum} to the free surface boundary conditions
\begin{align}\label{eq:boundary_condition_freeSurface}
 p\normal_\Sigma - \mu (\nabla \vel + \nabla \vel^{\sf T}) \normal_\Sigma &= \sigma \kappa \normal_\Sigma.
\end{align}
Such an approximation reduces the computational effort as we do not need to solve the integral equations \eqref{eq:integral_mass_conservation} and \eqref{eq:integral_momentum_conservation} for the pressure and velocity fields in the gas phase. While the interface mesh is moved, the bulk mesh has to follow these deformations in order to maintain acceptable bulk mesh quality. For this particular ALE approach, this is achieved by solving a diffusion equation for the mesh motion 
\begin{align}
 \nabla \cdot (D \nabla \boldsymbol{w}) &= 0   \quad \text{in}~ \domain \setminus \Sigma \\
 \boldsymbol{w} \cdot \normal_\Sigma = V^\Sigma \quad \text{and} \quad   \textbf{P}_\Sigma \boldsymbol{w} &= 0  \quad \text{on} ~\Sigma,~ 
  \boldsymbol{w} \cdot \normal_\Gamma = \vel \cdot \normal_\Gamma \quad  \text{on} \, \Gamma,
\end{align}
where $\boldsymbol{w}$ is the mesh velocity and  $\textbf{P}_{\Sigma} = \textbf{I} - \boldsymbol{n}_\Sigma \otimes \boldsymbol{n}_\Sigma$ is the interface projection tensor. A non-homogeneous anisotropic diffusion coefficient $D$ can be used to obtain a problem-specific mesh motion. This can be used to, e.g. obtain a less deformed mesh near the interface or the walls. The corresponding boundary conditions guarantee, that the mesh follows the interface by satisfying the kinematic condition \eqref{eq:kinematic_condition}. The boundary condition at the contact line ensures that the mesh representing the contact line remains on the boundary of the domain, e.g. a solid wall.

\begin{sloppypar}
The outlined ALE method has been implemented in OpenFOAM extend. The corresponding discretization of the free surface boundary conditions and details on the implementation can be found in \cite{Tukovic2012}. The curvature is computed by application of the surface divergence theorem to a face that is part of the discrete interface mesh  \cite{Tukovic2012}.
A finite volume method is used to discretize the mass, the momentum and the mesh motion equations. The pressure velocity coupling is solved by a segregated approach with a pressure Laplace equation. For the discretized pressure Laplace equation a \textit{Rhie-Chow-correction} is used \cite{Rhie1983}. This method corresponds to extending the stencil for the discretized Laplace operator aiming to eliminate oscillations in the discrete pressure field. The interface itself is moved by a \textit{control-point-algorithm} introduced in \cite{Muzaferija1997}. This approach attempts to minimize the phase volume conservation errors that typically occur in the numerical simulations of multiphase flows. For this approach, the stability condition $\Delta t < \sqrt{\rho \min(L_{PeN})^3/(2\pi \sigma)}$ has to be satisfied, where the discrete geodesic distance between cell centers $P$ and $N$ is denoted by $L_{PeN}$ \cite{Tukovic2012}. This condition is based on the criterion provided in \cite{Brackbill1992}.
The existing implementation has been extended to deal with wetting phenomena. This allows to employ a contact angle model that depends on the capillary number. The contact angle is adjusted by applying a suitable force at the contact line. This method adapts the curvature close to the contact line and thereby yields an interface that is intersecting the wall at the prescribed angle.
\end{sloppypar}

As wall boundary conditions, a combination of no penetration and Navier slip boundary condition has been implemented. The implementation reduces to the corresponding no penetration no slip condition for vanishing slip length. This general approach allows to use exactly the same implementation for interTrackFoam as for the interFoam solver described below.

    \subsection{Volume of Fluid method - interface capturing}
    
\paragraph{The VOF transport equation:} The Volume-of-Fluid method introduced in \cite{Hirt1981} uses the phase indicator function
\[ \chi(t,x) := \begin{cases} 1 & \text{if} \quad x \in \domain_l(t),\\ 0 & \text{if} \quad x \notin \domain_g(t) \end{cases}  \]
to track the fluid-fluid interface. With the help of the phase indicator function, the kinematic boundary condition \eqref{eq:kinematic_condition} can be expressed as the hyperbolic transport equation
\begin{align}
\label{eq:transport_phase_indicator}
\frac{D}{Dt} \, \chi(t,x) = 0 \quad \Leftrightarrow \quad \partial_t \chi(t,x) + \boldsymbol{u}(t,x) \cdot \nabla \chi(t,x) = 0.  
\end{align}
Physically, it states that the fluid particles cannot cross the fluid-fluid interface. Hence, in this model $\Sigma(t)$ is a material interface. But note that \eqref{eq:transport_phase_indicator} has to be understood in the sense of distributions since $\chi$ is discontinuous. This is a technical disadvantage compared to other methods such as the level set method, which use a smooth function to track the interface as its zero contour. Clearly, the level set function also satisfies the transport equation \eqref{eq:transport_phase_indicator}, where all derivatives can be evaluated pointwise. However, the advantage of the VOF method is that it automatically ensures conservation of the total volume in the discrete case: Since $\boldsymbol{u}$ is divergence free, it follows by integration of \eqref{eq:transport_phase_indicator} over a fixed control volume that
\[ \frac{d}{dt} \int_V \chi(t,x) dV = - \int_{\partial V} \chi(t,x) \boldsymbol{u}(t,x) \, dA. \]
Hence, every consistent numerical flux will conserve the total volume of the fluid.

\paragraph{The CSF model:} The ``Continuum Surface Force'' (CSF) model developed in  \cite{Brackbill1992} allows to treat the two-phase flow problem described in Section \ref{subsec:govEqs} as a single fluid with \emph{discontinuous} material properties $\rho$ and $\mu$ given by
\begin{align}
\label{eqn:discontinuous_material_properties}
\rho = \chi \rho_l + (1-\chi) \rho_g, \quad \mu = \chi \mu_l + (1-\chi) \mu_g. 
\end{align}
The effect of surface tension is modeled by a singular source term at the interface. The CSF model reads
\begin{align}
\rho (\partial_t \boldsymbol{u} + \boldsymbol{u} \cdot \nabla \boldsymbol{u}) - \mu \Delta \boldsymbol{u} + \nabla p = \boldsymbol{f} + \sigma \kappa \nsigma \delta_\Sigma \quad \text{in} \quad \domain,\label{eq:csf_1}\\
\nabla \cdot \boldsymbol{u} = 0 \quad \text{in} \quad \domain,\label{eq:csf_2}
\end{align}
where $\delta_\Sigma$ denotes the surface delta distribution. Hence, the full model in the CSF formulation is given by the equations \eqref{eq:transport_phase_indicator},\eqref{eq:csf_1},\eqref{eq:csf_2},\eqref{eq:navier_condition_v2},\eqref{eq:contact_angle_boundary_condition}. To satisfy the transmission condition \eqref{eq:jump_condition_velocity}, only globally continuous solutions for the velocity are allowed.

        \subsubsection{Geometric Volume of Fluid - FS3D}
The two-phase flow solver \emph{Free Surface 3D} (FS3D), originally developed by \cite{Rieber2004}, is based on a geometrical VOF method to discretize the incompressible two-phase Navier Stokes equations. Here, we give a brief overview of the main part of the discretization. More details on the methods and implementations can be found in \cite{Rieber2004,Fath2015,Fath2016}.

The governing equations in the CSF formulation are discretized using the finite volume approach on a fixed Cartesian grid, where the grids for velocity and pressure are staggered to enhance the stability of the method \cite{Harlow1965}. The material parameters $\rho$, $\mu$ are volume averaged according to \eqref{eqn:discontinuous_material_properties}, where the phase indicator function is replaced by the discrete volume fraction in each cell. The time integration is based on an explicit Euler method where the pressure-velocity coupling is realized by Chorin's projection method which leads to an elliptic equation for the pressure to perform a projection onto the space of solenoidal velocity fields.

The interface is locally reconstructed as a plane in each cell (also known as PLIC) \cite{Rider1998}, where the method by \cite{Youngs1984} is used to estimate the interface normal vector based on the volume fraction field. Using the reconstructed interface geometry, the numerical fluxes for the volume fraction are computed using an operator splitting method \cite{Strang1968}. The surface tension force in \eqref{eq:csf_1} is discretized with the balanced CSF method introduced in \cite{Popinet2009}. A height function representation of the interface is constructed in order to approximate the mean curvature. It has been demonstrated that this method is able to significantly reduce spurious currents at the interface.

Following the approach by \cite{Afkhami2008,Afkhami2009}, the height function is also used to indirectly enforce the contact angle boundary condition. The idea is to extrapolate the height function at the contact line linearly into a ghost cell layer, where the slope of the extrapolated interface is determined by the desired contact angle. As a result, the approximated value of the mean curvature is altered leading to a ``numerical force'', which drives the interface towards the desired contact angle. A drawback of that method is that it may create spurious currents at the contact line.

\begin{figure}[h]
 \centering
 \subfigure{\includegraphics{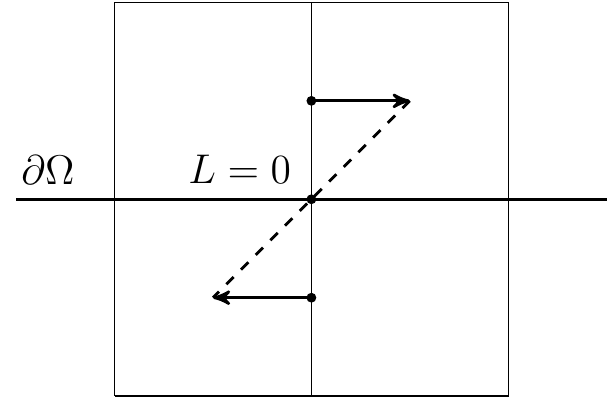}}
 \subfigure{\includegraphics{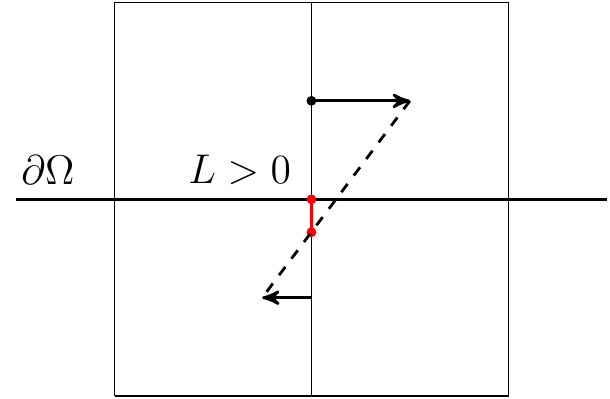}}
 \caption{Ghost cell based numerical realization of the no-slip and slip boundary conditions in FS3D.}
 \label{fig:slip_in_fs3d}
\end{figure}

To allow for a motion of the contact line one can either make use of the \emph{numerical slip} inherent to the method or prescribe the Navier slip condition \eqref{eq:navier_condition_v2}. The numerical slip is a property of the advection algorithm, which uses face-centered values to transport the volume fraction field. The face-centered velocity at the boundary cell layer is not identically zero and therefore allows for the motion of the contact line (see Figure~\ref{fig:slip_in_fs3d}). However, the numerical slip decreases with the mesh resolution leading to a significant mesh dependence of the solution. This numerical effect has been first described in the context of VOF methods in \cite{Renardy.2001}.

At the in- and and outflow boundary of the capillary, the pressure is fixed to zero and the components of the velocity field are copied into the dummy cell layers which formally corresponds to $\partial_n \vec{u} = 0$.
In order to ensure stability of the numerical method, the timestep is chosen according to the stability criterion
\[ \Delta t = \min \{ (\Delta t)_\sigma, (\Delta t)_\mu, (\Delta t)_{\boldsymbol{u}} \} \]
with the timescales (see \cite{Tryggvason2011} for a similar criterion) given by
\[ (\Delta t)_\sigma = \sqrt{\frac{(\rho_l+\rho_g)(\Delta x)^3}{4\pi\sigma}}, \quad (\Delta t)_\mu = \frac{\rho_l (\Delta x)^2}{6 \mu_l}, \quad (\Delta t)_{\boldsymbol{u}} = \frac{\Delta x}{\lVert \boldsymbol{u} \rVert}_\infty.  \]

        \subsubsection{Algebraic Volume of Fluid - interFoam}  
        The main difference between a geometric and an algebraic VOF approach is how the interface is advected. Since with the algebraic approach no reconstruction of the interface is used, an unstructured mesh discretization is straightforward. Hence, more complex geometries do not require any special treatment. Geometric VOF algorithms on arbitrary unstructured meshes are under active development \cite{Ivey2017,Jofre2015,Maric2018}. The algebraic VOF solver interFoam is implemented in OpenFOAM-based on unstructured meshes and is used in the present study.


In order to counteract numerical diffusion at the interface, an artificial compression term is added to the volume averaged form of the advection equation \eqref{eq:transport_phase_indicator}, leading to \cite{Weller2006}
\begin{equation}\label{eq:advectionCounterGrad}
  \frac{\partial \alpha }{\partial t } + \nabla\cdot(\alpha \boldsymbol{u}) + \nabla \cdot \left(\alpha \left(1-\alpha\right)\boldsymbol{u}_c\right)= 0\,.
\end{equation}
Here, $\alpha$ represents the volume-fraction field of phase one within each control volume. The compression velocity \(\boldsymbol{u}_c\) is given by \(\boldsymbol{u}_c = |\boldsymbol{u}| \frac{\nabla \alpha }{ |\nabla \alpha| }\). Equation \eqref{eq:advectionCounterGrad} is solved by the Multidimensional Universal Limiter with Explicit Solution (MULES), which preserves the boundedness of the volume fraction field \(\alpha\). The pressure-velocity coupling is solved using the Pressure Implicit Splitting Operator Algorithm (PISO) developed by \cite{Issa1986}.

\begin{figure}[H]
\centering
\includegraphics[width=.9\textwidth]{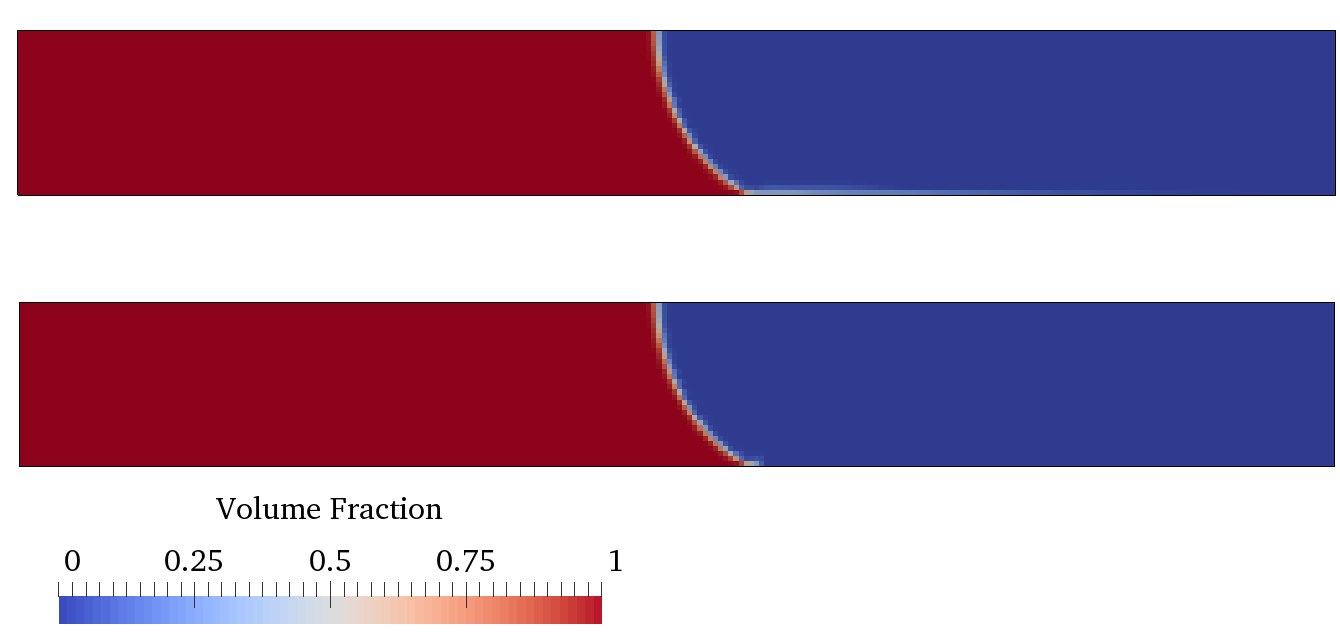}
\caption{Phase fraction for the standard contact angle treatment (top) in comparison with a localized contact angle model (bottom)}
\label{fig:interFoam:smearedCL}
\end{figure}

Before the curvature can be calculated, the interface normal at the wall boundary faces has to be corrected, to comply with the target contact angle \(\theta\), satisfying \(\boldsymbol{n}_\mathrm{\Sigma,w} \cdot \boldsymbol{n}_\mathrm{w} = \cos \theta\) at the contact line. Usually, this correction is applied along the entire wall boundary, even where the contact line is not present. However, in combination with the Navier slip model employed in this work, this procedure leads to excessive and unphysical smearing of the liquid-gas interface, as can be seen in Figure \ref{fig:interFoam:smearedCL} (top). Instead, using a heuristic approach, the contact angle is corrected in a band of approximately eight cells around the contact line. Outside this band the interface normals at the wall are corrected corresponding to a contact angle of $\SI{90}{\degree}$. Figure \ref{fig:interFoam:smearedCL} (bottom) shows the volume fraction field $\alpha$ for this correction method. By incorporating the corrected interface normal into the curvature calculation, the surface tension at the wall corrects the interface normal towards the target contact angle within the next time step.

On unstructured meshes, the use of  height functions to accurately compute the curvature of the interface \(\kappa\) is typically regarded inefficient. Therefore, two approaches to calculate the curvature are given here. First, the interface normal is approximated by the gradient of the volume fraction, $\boldsymbol{n}_\Sigma = \nabla \alpha / |\nabla \alpha|$. Then, the curvature is computed as  
\begin{equation}
    \kappa = - \nabla \cdot \boldsymbol{n}_\Sigma = - \nabla \cdot \frac{\nabla \alpha}{ |\nabla \alpha|} \,.
\end{equation}
\begin{sloppypar}
This approach provides a computationally cheap, but rather poor estimate of the interface normal \cite{Popinet2018}. Better accuracy can be achieved by reconstructing a volume fraction iso-surface \(\alpha = 0.5\) \cite{Kunkelmann2011,Batzdorf2015,Deising2018,Thammanna2018}. At each time step, the cell centered \(\alpha\)-values are interpolated to the cell-vertices. Linear interpolation of the interpolated volume fraction fields at the cell-vertices yields the points \(\boldsymbol{x}_p\), where the iso-contour cuts a cell-edges. Following the geometrical considerations of \cite{Lopez2008}, the cutting points \(\boldsymbol{x}_p\) are used to form an oriented surface element with the surface area \(|\boldsymbol{S}_\Sigma|\) and an interface normal \(\boldsymbol{n}_\Sigma = \boldsymbol{S}_\Sigma /|\boldsymbol{S}_\Sigma|\). This approach is used for the interFoam solver throughout the present work.
\end{sloppypar}
Finally, the surface tension is calculated following the CSF model as described by equation \eqref{eq:csf_1}, where $\delta_\Sigma$ is approximated by
\begin{equation}
\nsigma \delta_\Sigma \approx \nabla \alpha.
\end{equation}        
Thereby, the discretization scheme for the pressure gradient as well as for the of the volume fraction gradient are equal in order to maintain a `pseudo-well-balancedness' \cite{Yamamoto2016,Popinet2018}.
  
For the algebraic VOF simulations presented in the following sections, the time step was set according to
\begin{equation}
(\Delta t)_\mathrm{aVOF} = \mathrm{min} \{ \num{0.5}(\Delta t)_{\sigma,\rho_\mathrm{l}},\num{0.4}(\Delta t)_{\boldsymbol{u}} \}
\end{equation}
with
\begin{equation}
(\Delta t)_{\sigma,\rho_\mathrm{\ell}}=\sqrt{\frac{\rho_\mathrm{\ell}(\Delta x)^3}{4\pi\sigma}}<\sqrt{\frac{(\rho_\mathrm{l}+\rho_\mathrm{g})(\Delta x)^3}{4\pi\sigma}},
\end{equation}
hence satisfying the criterion given by \cite{Brackbill1992}.

    \subsection{Extended Discontinuous Galerkin - BoSSS}
    The code-framework BoSSS (\textit{Bounded Support Spectral Solver}), is based on the Extended Discontinuous Galerkin method (Extended DG/XDG) \cite{kummer_bosss_2012}, also referred to as Unfitted DG \cite{heimann_unfitted_2013}
or Cut-Cell DG method. 
The interface is represented by a level-set function $\varphi(\boldsymbol{x},t)$, such that
\begin{equation}
	\Sigma(t) = \{ \boldsymbol{x} \in \mathcal{D}: \varphi(\boldsymbol{x}, t) = 0 \}, \quad
	\domain_l(t) = \{ \boldsymbol{x} \in \mathcal{D}: \varphi(\boldsymbol{x}, t) < 0 \}, \quad 
	\domain_g(t) = \{ \boldsymbol{x} \in \mathcal{D}: \varphi(\boldsymbol{x}, t) > 0 \}.
\end{equation}
The starting point of the XDG discretization is a standard DG polynomial basis 
$\boldsymbol{\Phi} = (\phi_{j,n})_{j = 1,\ldots,J \ n=0,\ldots, N_p}$.
Here, $J$ and $N_p$ denote the number of cells and the number of polynomials up to degree $p$, respectively.
The individual functions $\phi_{j,n}(\boldsymbol{x})$ are defined as polynomials 
within the computational cell $K_j$ and zero outside.
In the XDG approach  \cite{kummer_extended_2016}, the DG approximation space is adapted in every time-step in order to be conformal with the fluid interface. This allows e.g. a sharp representation of the pressure jump.
The adaptation can be obtained by multiplying the basis functions $\phi_{j,n}$ with the Heaviside 
function of $\varphi$, i.e. some field quantity $q(\boldsymbol{x})$
is approximated as a weighted sum 
\begin{equation}
q(\boldsymbol{x},t) = 
 \sum_{j,n} 
       q_{j,n}^-(t) \phi_{j,n}(\boldsymbol{x}) H(-\varphi(\boldsymbol{x}, t)) 
     + q_{j,n}^+(t) \phi_{j,n}(\boldsymbol{x}) H(+\varphi(\boldsymbol{x}, t)) .
\end{equation}  
In each cut-cell (a cell intersected by the interface $\Sigma$) this yields
separate degrees of freedom $q_{j,n}^-(t)$ and $q_{j,n}^+(t)$ for the  liquid and gaseous domain, $\domain_l$ and $\domain_g$, respectively. 

The level-set field $\varphi$ is represented by a standard DG approach, i.e.
\begin{equation}
 \varphi(\boldsymbol{x},t) = 
\sum_{j,n} \varphi_{j,n}(t) \phi_{j,n}(\boldsymbol{x}) .
\end{equation}
Its evolution is implemented as a fast-marching procedure in a narrow band around the interface. With this approach an extension-velocity problem is solved \cite{utz_high-order_2018} and if necessary a reinitialization \cite{utz_interface-preserving_2017} is performed in each cell.

Since the XDG basis functions 
$\phi_{j,n}(\boldsymbol{x}) H( \pm \varphi(\boldsymbol{x}, t))$ are conformal with the fluid interface,
XDG allows a sharp and sub-cell accurate representation of kinks in the velocity field as well as discontinuities in the pressure. In combination with the hierarchical-moment-fitting technique for the numerical integration on cut-cells \cite{muller_highly_2013}, this ansatz regains spectral convergence properties even for such low-regularity solutions. 

For the spatial discretization, we are considering the momentum balance equations in divergence form, cf. \eqref{eq:momCons}.
The convective terms are discretized by a local Lax-Friedrichs flux, as used in the work of \cite{shahbazi_high-order_2007}. For the viscous terms a variant of the standard symmetric interior penalty (SIP), first proposed by \cite{arnold_interior_1982}, is employed. In order to comply with the Lady\u{z}enskaja–Babu\u{s}ka–Brezzi condition (e.g. \cite{babuska_finite_1973,brezzi_existence_1974}), the velocities are discretized with order $k$ and the pressure field with $k^\prime = k - 1$. A complete description is given in \cite{kummer_extended_2016}. The discretization in time is done with an explicit moving interface approach following the work of \cite{kummer_time_2018}, where a BDF2 scheme is used.
 
Regarding the numerical computation of the surface tension force, a Laplace-Beltrami formulation 
\begin{equation}
    \int_{\Sigma} \text{div}_{\Sigma}\left( \sigma \textbf{P}_{\Sigma}\right) \cdot \boldsymbol{v} \text{d}S = - \int_{\Sigma} \sigma \textbf{P}_{\Sigma} : \nabla_{\Sigma} \boldsymbol{v} \text{d}S + \int_{\partial \Sigma} \sigma \textbf{P}_{\Sigma} \boldsymbol{\tau}_{\partial \Sigma} \cdot \boldsymbol{v} \text{d}l,
\label{eq:stf_LB}	
\end{equation}
is used, where $\nabla_{\Sigma} = \textbf{P}_{\Sigma} \nabla$. Concerning dynamic contact line problems, this approach allows the introduction of the generalized Navier boundary condition (GNBC) into the spatial discretization \cite{reusken_finite_2017}.
The GNBC extends the Navier boundary condition \eqref{eq:navier_condition} with a force balance at the contact line $\partial \Sigma$:
\begin{equation}
	\sigma \left( \cos{\theta} - \cos{\theta_\text{e}} \right) \boldsymbol{n}_{\cl} = -\fcL \left( \boldsymbol{u} \cdot \boldsymbol{n}_{\cl} \right) \boldsymbol{n}_{\cl}
	\label{eq:GNBC}
\end{equation}
where $\theta_\text{e}$ denotes the static equilibrium angle derived from Young's relation. Setting $\fcL = 0$ this boundary condition results in a static contact angle model, which is used throughout this paper. Furthermore the time step is always set according to the criterion provided by \cite{Brackbill1992} with $\Delta t < \sqrt{(\rho_l+\rho_g)(\Delta x)^3/4 \pi \sigma}$, where the DG-mesh size is defined by ${\Delta x}_{DG} = \Delta x / (k+1)$.

\section{Numerical setups and results} \label{sec:setups_results}
The continuum mechanical problem described in Section~\ref{subsec:govEqs} has been solved using the four numerical methods described in Section~\ref{sec:numerical_methods}. Using such diverse numerical approaches requires specific preparation of the corresponding initial conditions. These are described in the following subsections. All implementations use a dimensional formulation with the parameters given in Section~\ref{sub:validation_case}. The simulation results are compared with the ODE model described in Section~\ref{sub:ode_model} using the scaling introduced in Section~\ref{sub:scaling}.

    \subsection{Physical parameters}\label{sub:validation_case}
    We aim to compare the results of the four numerical approaches for the continuum mechanical problem with the ODE model given in Table \ref{tab:scaledModels}. Various restrictions apply to the choice of physical parameters: 
\begin{enumerate}[(i)]
 \item The interface has to maintain its circular shape during the rise. This is an essential assumption for the derivation of the ODE model.
 \item The setup must yield a rise height that differs significantly from its initial height. This way, various rise behaviors known from literature can form.
 \item The simulation domain should be small to reduce computational costs.
 \item The physical parameters should be chosen such that the necessary number of time steps is as small as possible.
 \item The influence of the gas phase should be small in order to ensure compatibility to the reference model.
\end{enumerate}
Following \cite{Fries2009}, we choose a set of physical parameters (see Table~\ref{tab:physicalParams}) to vary the non-dimensional group 
\[ \Omega = \sqrt\frac{9 \sigma \cos \theta \mu^2}{\rho^3 g^2 R^5} \]
in the range of $\Omega = 0.1, 0.5, 1, 10, 100$ in order to cover all rise regimes from a highly oscillatory to a strictly monotone rise. We fix the contact angle for all cases to $\theta=30^\circ$. Moreover, it is convenient to fix the computational domain. We, therefore, enforce the same stationary height for all cases via
\begin{align}
\label{eqn:eotvos_relation}
\hjur = 4 R \quad \Leftrightarrow \quad 4 = \frac{\hjur}{R} = \frac{\cos\theta}{\text{Eo}}. 
\end{align}
The initial height is fixed to $h_0 = 2R = \hjur/2$. Since the rise height is expected to show strong oscillates for small $\Omega$ numbers, the height of the computational domain is set to $h_\text{D} = 2 \hjur =  8R$ for all implementations following a two phase approach. For the two-phase flow solvers, the density and viscosity ratios are set to $\frac{\rho}{\rho_g} = 1000$ and $\frac{\mu}{\mu_g} = 1000$ to ensure that the influence of the gas phase is small.

In order to fulfill the first requirement in the stationary state (see \cite{Concus1968}), the half gap width $R$ should be small compared to the capillary length, i.e.
\begin{align}\label{eq:condRadius}
 \frac{R}{\lcap} = \sqrt{\frac{\rho g R^2}{\sigma}}  = \sqrt{\text{Eo}} 
\end{align}
should be small. According to \eqref{eqn:eotvos_relation}, the Eotvos number is fixed to $\text{Eo} = \cos(30^\circ)/4 \approx 0.217$. Moreover, the shape of the meniscus can be expected to remain circular during the rise if the Capillary number
\[ \text{Ca} = \frac{\mu V_{cl}}{\sigma} \]
is small enough. This quantity cannot be computed directly since the contact line velocity $V_{cl}$ is not known a priori. However, the maximum capillary numbers extracted from the performed simulation show that this requirement is fulfilled for the considered parameters, see Table~\ref{tab:physicalParams}.

\begin{table}[H] 
	\caption{Physical parameters for the $\Omega$-study.}
	\label{tab:physicalParams}
	\centering
	\begin{tabular}{ l l l l l l l l l}
		$\Omega$ & $R $ & $\rho $ & $\mu $ & $g$ & $\sigma$ & $\theta_{\text{e}}$ &  $Ca_\text{max}$ & Eo \\ 
		~- & \SI{}{\meter} & \SI{}{\kilogram \per \meter \cubed} & \SI{}{\pascal \second} & \SI{}{\meter \per \second \squared} & \SI{}{\newton \per \meter} & \SI{}{\degree} & ~- & ~- \\
		\hline
		0.1 & 0.005 & 1663.8 & 0.01 & 1.04 & 0.2 & $30$ & 0.003 & 0.217 \\ 
		0.5 & 0.005 & 133.0 & 0.01 & 6.51 & 0.1 & $30$ & 0.015 & 0.217 \\ 
		1 & 0.005 & 83.1 & 0.01 & 4.17 & 0.04 & $30$ & 0.029 & 0.217 \\ 
		10 & 0.005 & 3.3255 & 0.01 & 26.042 & 0.01 & $30$ & 0.106 & 0.217 \\ 
		100 & 0.005 & 0.33255 & 0.01 & 26.042 & 0.001 & $30$ & 0.110 & 0.217
	\end{tabular}
\end{table}

\paragraph{Viscous vs. capillary timestep limit:}
To estimate the computational costs, we consider the numerical timestep limits already discussed in Section~\ref{sec:numerical_methods}. The timestep limit due to surface tension is given by\footnote{Note that the applied $(\Delta t)_\sigma$ might differ by a constant factor between the considered numerical methods, cf. Section~\ref{sec:numerical_methods}. Here we consider the actual values used for FS3D.}
$ (\Delta t)_\sigma = \sqrt{\rho (\Delta x)^3/(4\pi\sigma)}$.
Here $\Delta x$ denotes the mesh size which is assumed to be constant for the purpose of estimating the computational costs. Since the viscous terms are explicitly discretized in FS3D, there is also a viscous timestep limit for the latter method given by
$(\Delta t)_\mu = \rho (\Delta x)^2/(6\mu)$.
In order to identify which one is the limiting value for a given set of parameters, we take the ratio of the two quantities to get
\begin{align} 
\frac{(\Delta t)_\sigma}{(\Delta t)_\mu} = \frac{3}{\sqrt\pi} \frac{\mu}{\sqrt{\sigma \rho R}} \left(\frac{R}{\Delta x}\right)^{1/2} = \frac{3}{\sqrt\pi} \, \text{Oh} \, N_{cells}^{\frac12}, 
\end{align}
where $N_{cells} = R/\Delta x$ is the number of computational cells per radius and 
\[ \text{Oh} := \frac{\mu}{\sqrt{\sigma \rho R}} = \Omega \frac{\text{Eo}}{3 \sqrt{\cos\theta}} \] 
is the Ohnesorge number. Hence, the two timestep limits are equal for
\[ N_{cells}^\ast = \frac{\pi}{9 \, \text{Oh}^2} = \frac{\pi \cos \theta}{\Omega^2 \text{Eo}^2}. \]
The capillary timestep limit $(\Delta x)_\sigma$ is dominant for $N_{cells} \leq N_{cells}^\ast$, while the viscous timestep limit $(\Delta x)_\mu$ is dominant for $N_{cells} \geq N_{cells}^\ast$. Hence the viscous timestep limit dominates over the capillary timestep limit for large $\Omega$.

\paragraph{Capillary timestep limit in different regimes:}
Following \cite{Fries2009}, we consider the time-scales
\begin{align*}
 t_{\text{scale},1} = \frac{\rho^2 g^2 R^3}{3\mu \sigma \cos \theta}, \quad t_{\text{scale},2} = \sqrt{\frac{\rho g^2 R}{\sigma \cos \theta}}, \quad t_{\text{scale},3} = \frac{3\mu}{\rho R^2}
\end{align*}
Note that the dimensionless time $t^*$ is defined as $t^* = t_{\text{scale},k} \, t$ where $t$ is the physical time. We are interested in the number of timesteps necessary to compute up to dimensionless time $t^\ast=1$ in scaling $k$ which can be estimated either by
\[ N_{steps}^{k,\sigma} = \frac{1}{t_{\text{scale},k} (\Delta t)_\sigma}, \quad \text{or} \quad N_{steps}^{k,\mu} = \frac{1}{t_{\text{scale},k} (\Delta t)_\mu}. \]
A straightforward calculation shows the relations
\begin{align*}
N_{steps}^{1,\sigma} &\propto \, \frac{\cos \theta}{\text{Eo}} \, \Omega \, N_{cells}^{3/2}, \quad N_{steps}^{2,\sigma} \propto \frac{\sqrt{\cos\theta}}{\text{Eo}} N_{cells}^{3/2}, \quad N_{steps}^{3,\sigma} \propto \frac{\sqrt{\cos\theta}}{\text{Eo}} \, \frac{1}{\Omega} \, N_{cells}^{3/2}, \\
N_{steps}^{1,\mu} &\propto \sqrt{\cos\theta} \, \Omega^2 N_{cells}^2, \quad N_{steps}^{2,\mu} \propto \Omega \, N_{cells}^2, \quad N_{steps}^{3,\mu} \propto N_{cells}^2.
\end{align*}
Hence both large (for scaling I) and small (for scaling III) values of $\Omega$ lead to a large number of timesteps. In particular, the number of timesteps in scaling 1 grows like $\Omega^2$ for FS3D where the viscous timestep limit is essential for large $\Omega$. Table~\ref{tab:timesteps} gives an overview of the timesteps necessary for a uniform mesh with $32$ cells per radius.

\begin{table}[h]
\caption{Estimate of the number of timesteps for $N_{cells}=32$ and the physical parameters defined in Table~\ref{tab:physicalParams}.}
\centering
\begin{tabular}{l l l l l l l l}
$\Omega$ &  $N_{cells}^\ast$ & $N_{steps}^{1,\sigma}$ & $N_{steps}^{2,\sigma}$ & $N_{steps}^{3,\sigma}$ & $N_{steps}^{1,\mu}$ & $N_{steps}^{2,\mu}$ & $N_{steps}^{3,\mu}$ \\
\hline
0.1 &  5804 & $2.57 \cdot 10^2$  & $2.76 \cdot 10^3$ & $\mathbf{4.78 \cdot 10^4}$ & $1.91 \cdot 10^1$ & $2.05 \cdot 10^2$ & $3.55 \cdot 10^3$ \\
0.5 &  232 & $1.28 \cdot 10^3$ & $2.76 \cdot 10^3$ & $\mathbf{9.55 \cdot 10^3}$ & $4.76 \cdot 10^2$ & $1.02 \cdot 10^3$ & $3.55 \cdot 10^3$ \\
1 &  58 & $2.57 \cdot 10^3$ & $2.76 \cdot 10^3$ & $\mathbf{4.78 \cdot 10^3}$ & $1.91 \cdot 10^3$ & $2.05 \cdot 10^3$ & $3.55 \cdot 10^3$ \\
10 &  0.58 & $\mathbf{2.57 \cdot 10^4}$ & $2.76 \cdot 10^3$ & $4.78 \cdot 10^2$ & $\mathbf{1.91 \cdot 10^5}$ & $2.05 \cdot 10^4$ & $3.55 \cdot 10^3$ \\
100 & 0.006 & $\mathbf{2.57 \cdot 10^5}$ & $2.76 \cdot 10^3$ & $4.78 \cdot 10^1$ & $\mathbf{1.91 \cdot 10^7}$ & $2.05 \cdot 10^5$ & $3.55 \cdot 10^3$ \\
\end{tabular}
\label{tab:timesteps}
\end{table}

   \subsection{Arbitrary Lagrangian Eulerian method - interface tracking} \label{sub:res_ALE}

\begin{figure}[H]
\centering
\subfigure{\includegraphics[trim=0 0 0 -0.88cm,width=0.3\linewidth]{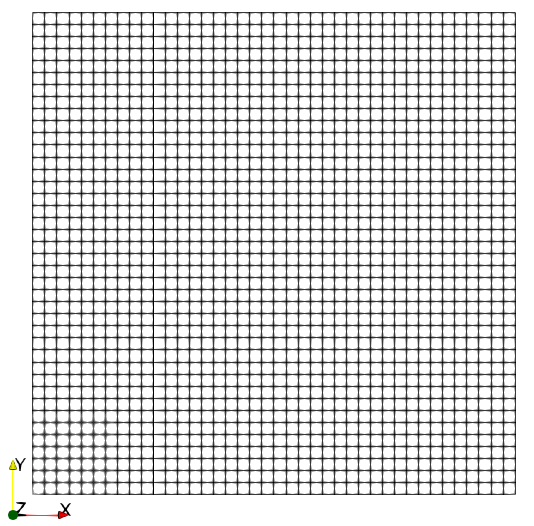}}
\subfigure{\includegraphics[width=0.3\linewidth]{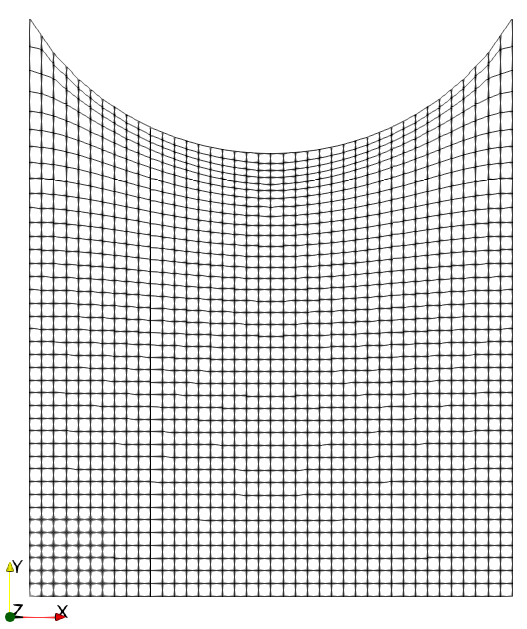}}
\caption{Left: Regular mesh with $40\times 40$ cells. Right: mesh used as initial condition.}
\label{fig:ale_initCond} 
\end{figure}

\paragraph{Numerical setup:} As the reference model assumes that the interface has the shape of a circular section, this initial condition has to be met by the initial condition of the interface tracking method. This means, that a mesh has to be provided such that the interface consists of mesh faces in the shape of a circular section. In addition, the mesh inside the domain has to provide sufficient quality such that the relevant fields can be resolved. The initial mesh is obtained as follows: First a square is tessellated with a cartesian mesh. Then, boundary conditions are applied as described in Subsection \ref{subsec:govEqs}, except for the inflow. Here, the inflow is first ``closed'' by choosing $\vel=0$ together with $\partial_n p=0$ on $\partial \mathcal{D}_{\text{in}}$. Then the simulation is started and run until an equilibrium state is obtained. This yields a reorientation of the interface to its static shape as shown in figure \ref{fig:ale_initCond}.

In order to resolve the velocity and pressure fields close to the interface, an anisotropic mesh diffusion coefficient is chosen. This allows to keep the mesh in interface vicinity more ``rigid'' and hence providing sufficient resolution near the contact line during the rise.

\begin{figure}[H]
\centering
\subfigure{\includegraphics[width=0.48\textwidth]{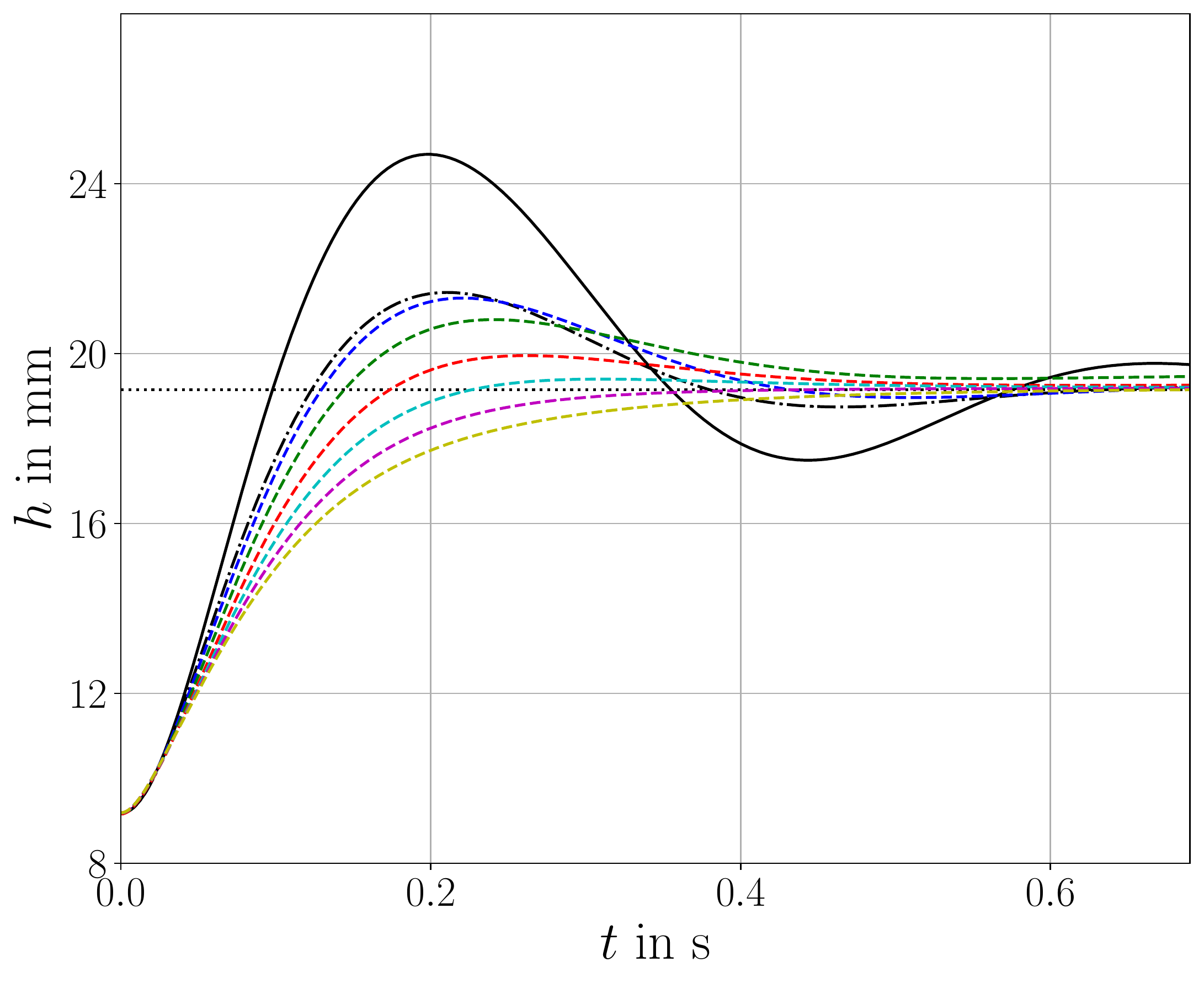}\label{fig:ale_meshConv_numSlip}} 
\subfigure{\includegraphics[width=0.48\textwidth]{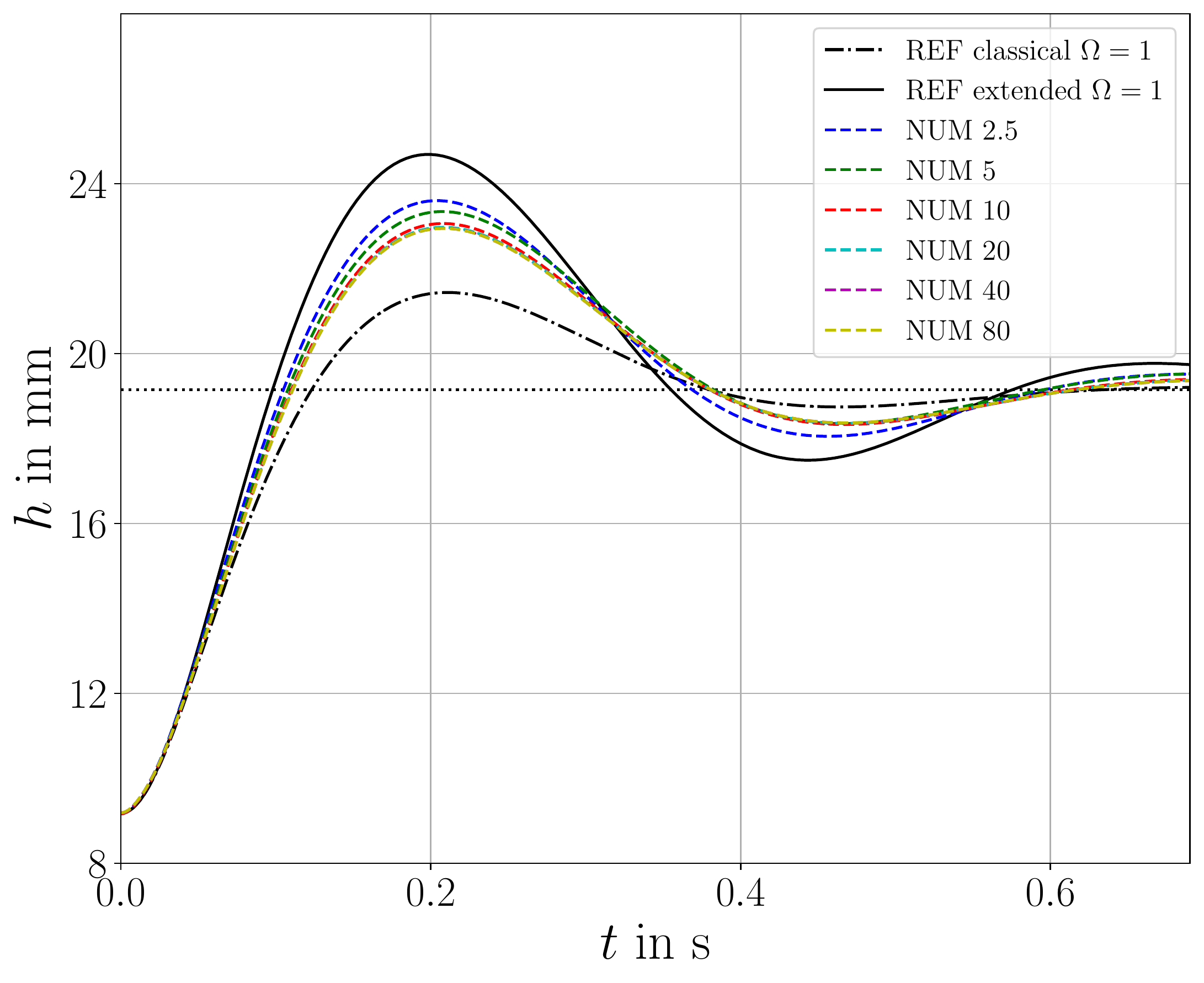}\label{fig:ale_meshConv_navSlip}}
\caption{Mesh convergence study for $\Omega=1$ and meshes with 2.5 to 80 cells per radius. Left: no slip boundary conditions discretized with ``numerical slip''. Right: Navier slip with $L=R/5$.}
\label{fig:ale_meshConvergence} 
\end{figure}

\paragraph{Convergence Study:} Figure \ref{fig:ale_meshConvergence} shows a mesh convergence study for the ALE interface tracking method with parameters yielding the case $\Omega=1$. Meshes with a resolution of 5 to 160 cells per diameter have been used. The results in Figure \ref{fig:ale_meshConv_numSlip} are obtained with no-slip boundary conditions where the interface is moved using ``numerical slip``. In the ALE context we use the term ''numerical slip`` for an approach where the velocity from a mesh cell is used to move the mesh point representing the contact line mesh. The ALE-solution using numerical slip with a mesh resolution of 5 cells per diameter agrees well with the classical model for the capillary rise problem that that is based on a no-slip boundary condition on the wall. With increasing resolution however, the deviations from the reference solution increase. While for 5 cells per diameter a small oscillation can be seen, the solution with a mesh resolution above 40 cells per diameter does not show this feature. Its resemblance is closer to monotonic behavior with large $\Omega$, then an oscillatory one. Hence, the change in the solution is not only a quantitative but also a qualitative change. 
In Figure \ref{fig:ale_meshConv_navSlip} the results for the same convergence study as before is shown except for a Navier slip boundary condition with $L=R/5$ is being applied at the capillary walls. 

The differences between the curves from the various meshes are smaller then those from the approach using numerical slip. Furthermore, the curves with a resolution of more than 20 cells are visually coinciding, indicating a mesh-converged solution. As a significant slip length has been used, it can be expected that the converged solution does not fully agree with the classical model based on a no slip assumption (REF classical). The converged results are located between the classical solution and the extended solution. This reference curve incorporates the influence of significant slip on the capillary walls shows in this case a less damped oscillation that the full continuum mechanical solution. Altogether, the no slip boundary condition - discretized with numerical slip - does not yield mesh convergent results for this case, while the same setup gives converged results on comparably low resolved meshes when a Navier slip boundary condition with a large slip length is used.

One might argue that at some point, the solution will converge as soon as the mesh is sufficiently refined. With no-slip boundary condition however, it is known that the local velocity field (and its dissipation) at the contact line is non-integrable singular \cite{Huh1971}. Hence we do not expect a mesh converged solution for no-slip boundary conditions that have been discretized using ``numerical slip''. The situation changes when a Navier slip boundary condition is used. The pressure as well as the viscous dissipation are at least integrable \cite{Huh1971}.   
        
    \subsection{Volume of Fluid method - interface capturing}

        \subsubsection{Geometric Volume of Fluid - FS3D}
\begin{figure}[H]
 \includegraphics[width=\textwidth]{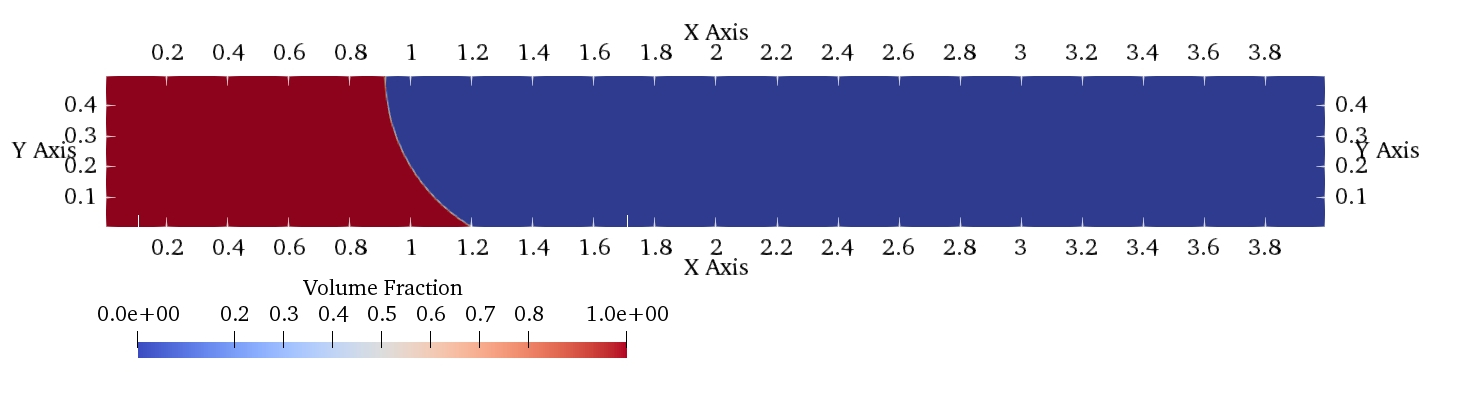}
 \caption{Simulation setup for FS3D. Length in \SI{}{\centi \meter}.}
 \label{fig:fs3d_setup}
\end{figure}

\begin{figure}[H]
\centering
\subfigure{\includegraphics[width=0.48\textwidth]{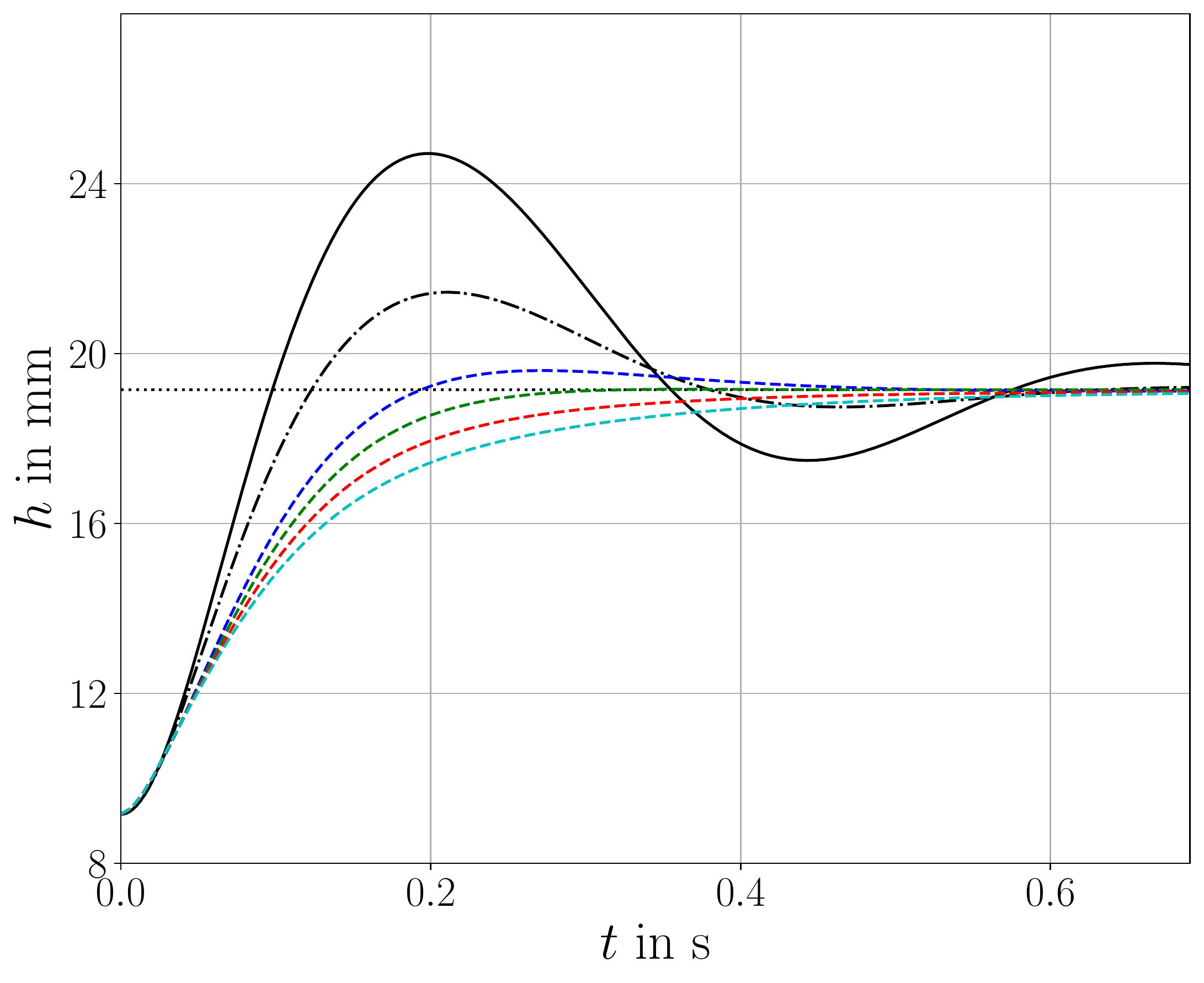}}
\subfigure{\includegraphics[width=0.48\textwidth]{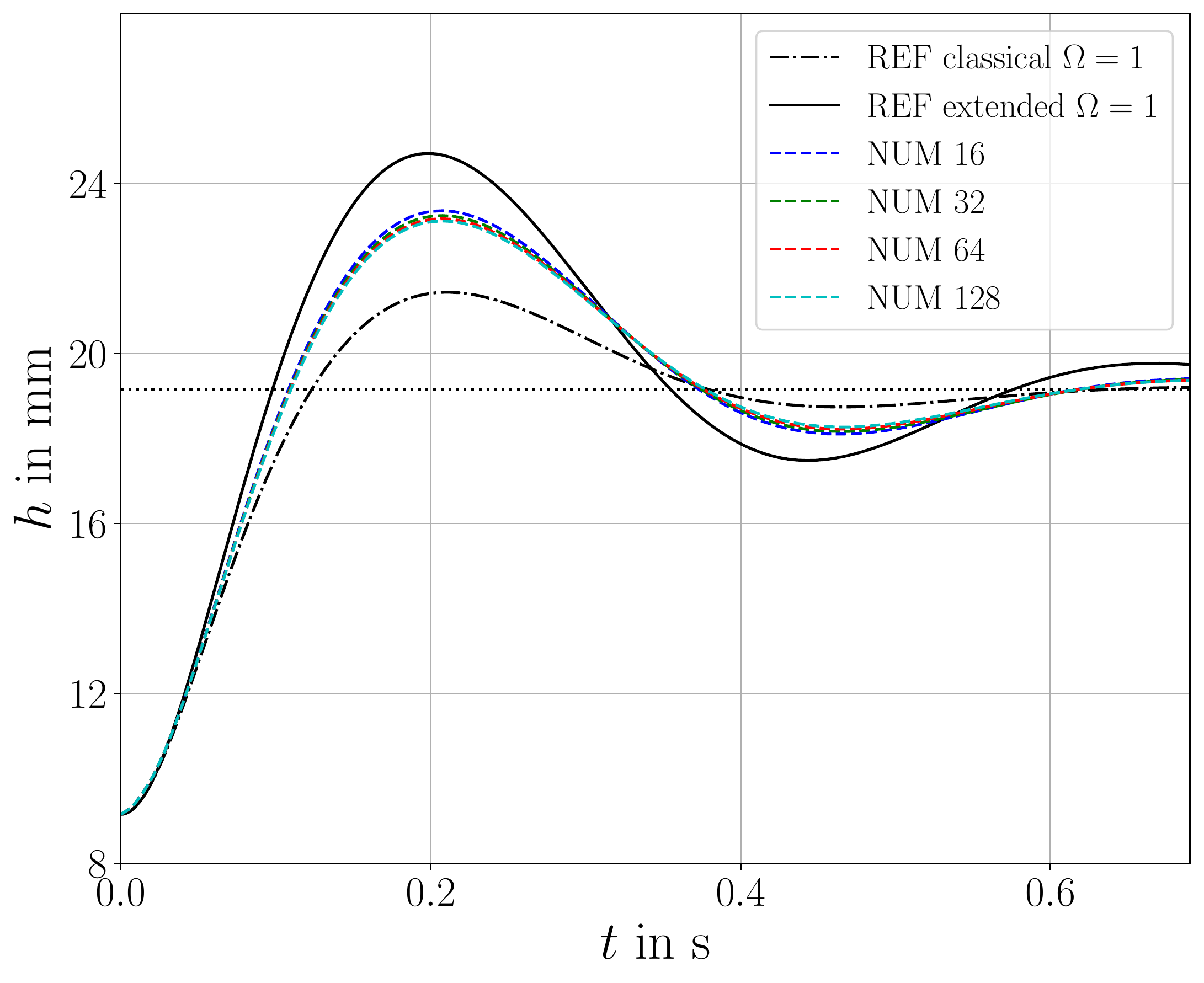}}
\caption{Mesh study with FS3D using numerical slip (left) and resolved Navier slip with $L=R/5$ (right).}
\label{fig:mesh_study_fs3d}
\end{figure}

\paragraph{Numerical setup:} We make use of the mirror symmetry of the problem and simulate only half of the domain on a uniform Cartesian grid with a symmetry plane in the center of the capillary, see Figure~\ref{fig:fs3d_setup}. The symmetry plane is realized as a second solid wall with 90 degree contact angle and free slip for the velocity, i.e. Equation~\eqref{eq:navier_condition} with $\lambda=0$. The volume fraction field is initialized for a spherical interface matching the contact angle boundary condition and a prescribed initial volume. The initial value for the velocity field is zero everywhere in the computational domain. The simulations in this study require between $4$ and $192$ cells in radial direction yielding $128$ to $294912$ cells in total.
\paragraph{Convergence study:} The results for the dynamics of the capillary rise are highly dependent on the mesh size if \emph{numerical slip} is used, see Figure~\ref{fig:mesh_study_fs3d}. This is due to the fact that the effective slip length of the method is related to the mesh size. In particular, with increasing mesh resolution the dynamics is increasingly damped. However, the final rise height is reasonably mesh independent and agrees well with the corrected Jurins Height given by \eqref{eq:stationary_rise_height} even if numerical slip is used. On the other hand, the dynamics of the apex height becomes mesh convergent if the Navier slip condition is applied with a macroscopic slip length $L=R/5$ which can be resolved by the computational mesh, see Figure~\ref{fig:mesh_study_fs3d}.

        \subsubsection{Algebraic Volume of Fluid - interFoam}
\begin{figure}[H]
\centering
\includegraphics[width=\textwidth]{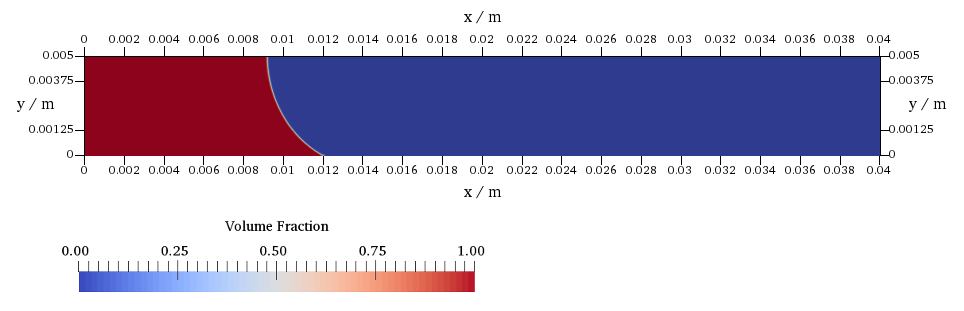}
\caption{Simulation setup for interFoam}
\label{fig:setup_interFoam}
\end{figure}

\paragraph{Numerical setup:} Similar to the setup for the FS3D solver, we also made use of the mirror symmetry for the simulations with interFoam. Accordingly, the same boundary conditions to those given above were used. Instead of imposing a contact angle of $\ang{90}$, however, a homogeneous Neumann boundary condition is set for the volume fraction field. Due to the no penetration condition at the symmetry plane, this is not relevant for the advection of the phase fraction. The boundary values of the phase fraction are evaluated for the interface normal and curvature calculation. Hence, the homogeneous Neumann condition ensures that the interface is perpendicular to the mirror plane. The initial shape of the meniscus, see Figure \ref{fig:setup_interFoam}, is found by initially simulating the capillary with a "closed" inlet, starting with an initially flat meniscus similar to the one displayed in Figure \ref{fig:ale_initCond} analog to the approach used for the ALE method. After reaching a stationary state the boundary conditions are changed to those described in Section \ref{subsec:govEqs}.

\begin{figure}[H]
\centering
\subfigure{\includegraphics[width=0.48\textwidth]{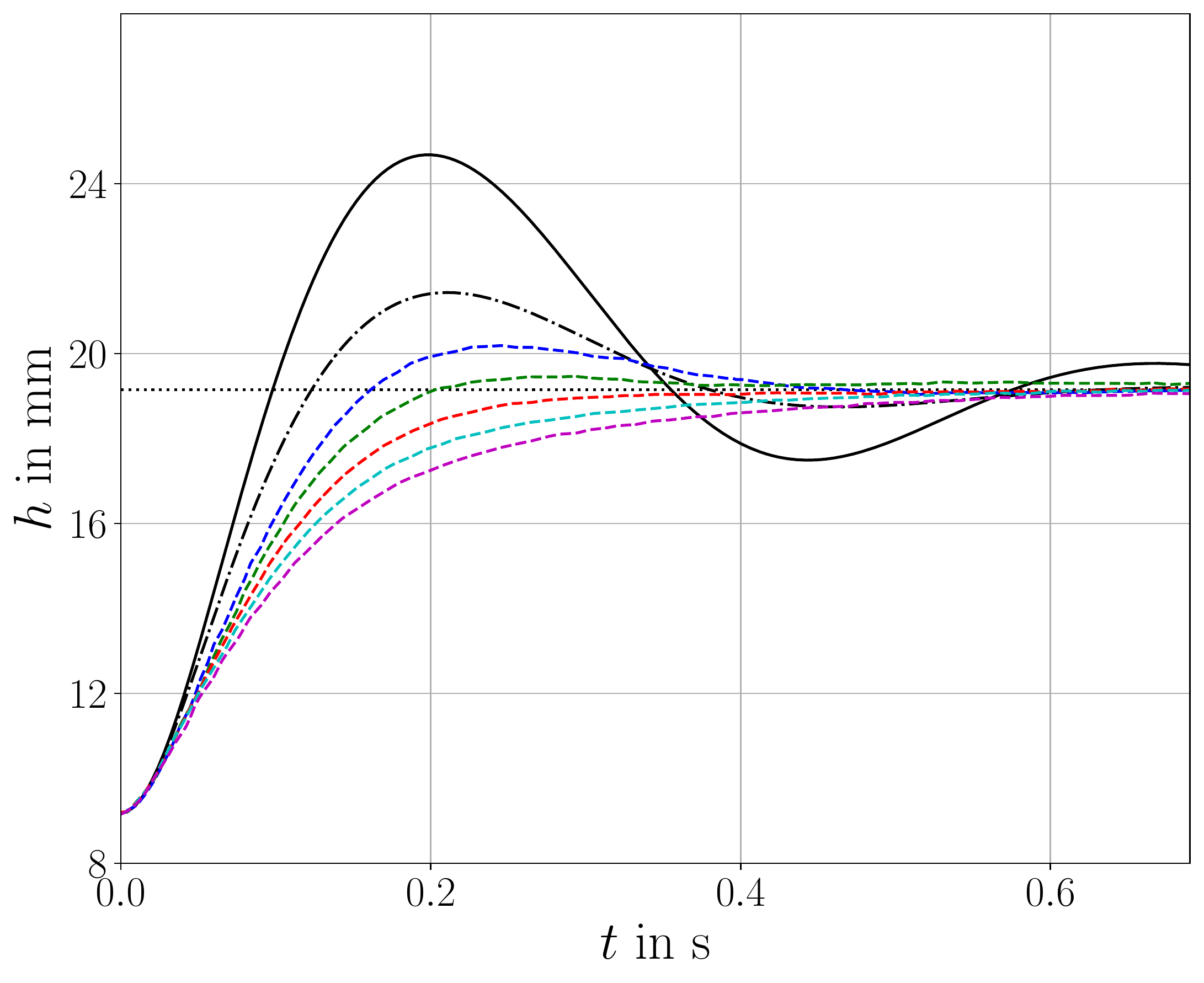}}
\subfigure{\includegraphics[width=0.48\textwidth]{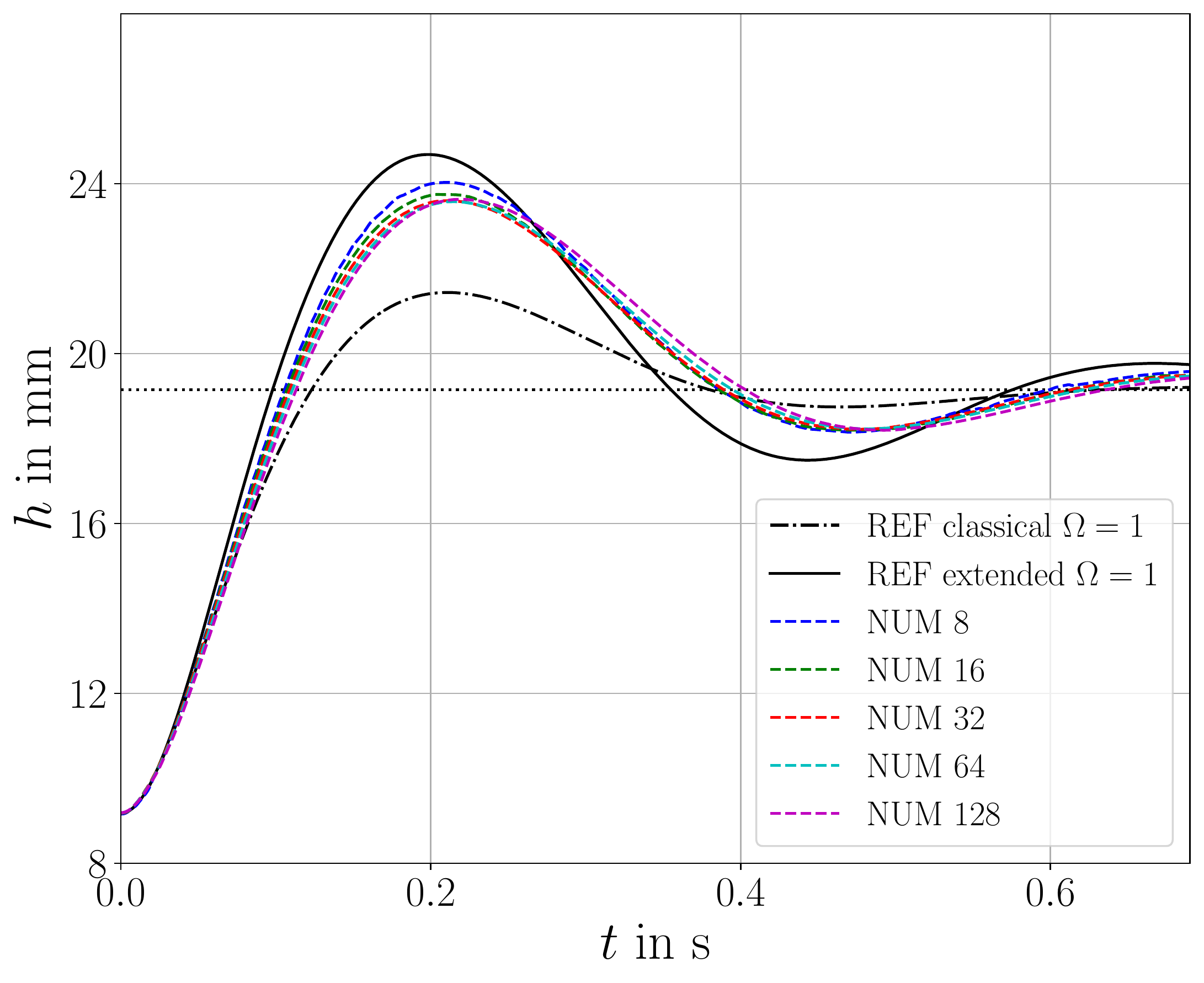}}
\caption{Mesh study with interFoam using numerical slip and Navier slip in comparison.}
\label{fig:MeshStudyInterFoam}
\end{figure}

\paragraph{Convergence study:} As for the other numerical methods, a mesh dependence study was conducted for the $\Omega=1$ case. Starting from a relatively coarse resolution of 8 cells per radius, the mesh size was halved up to a resolution of 128 cells per radius, i.\,e. the total number of cells ranging from 512 to 131072. Again, the mesh dependence was investigated for both numerical slip, as well as with the Navier slip condition at the solid wall. The results for the interFoam solver are displayed in Figure \ref{fig:MeshStudyInterFoam} in comparison with the solution of the 1D reference models. Similar to the results presented above, the rise dynamics are highly dependent on the mesh resolution if no slip boundary condition is used, i.\,e., if the effective slip length of the method is proportional to the mesh size (Figure \ref{fig:MeshStudyInterFoam} (left)). This mesh dependence is reduced, when the Navier slip model with resolved slip length is employed, as can be seen from Figure \ref{fig:MeshStudyInterFoam} (right). However, we did not find mesh convergence in this case. This might be due to increasing parasitic/spurious currents with an increasing mesh resolution - a well-known issue for VOF methods. Nevertheless, good agreement with the corrected Jurin's height $h_\infty^{apex}$ as given in equation \eqref{eq:stationary_rise_height} is found independent for both, the boundary condition and all spatial discretizations.

    \subsection{Extended Discontinuous Galerkin - BoSSS}
\begin{figure}[H]
\centering
\includegraphics[width=16cm]{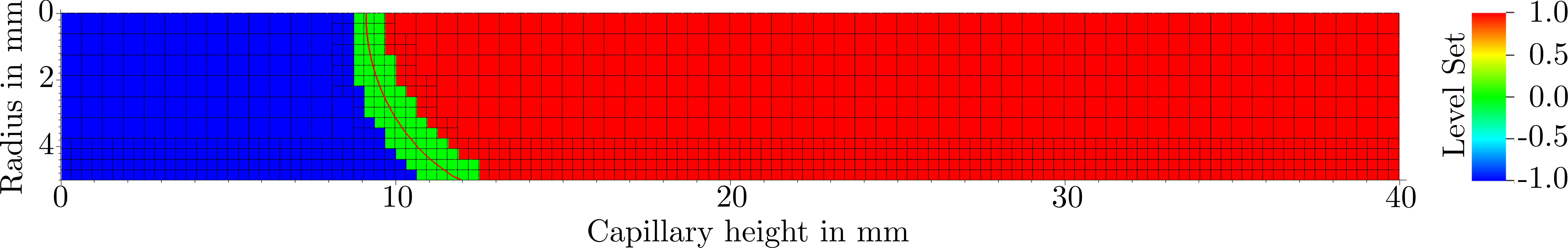}
\caption{Computational domain and initial condition after start-up phase. The level-set field $\phi$ is only computed in a narrow band around the interface (green cells). Outside the narrow band $\phi$ is set constant. In case of adaptive mesh refinement (here on level 1) the cells in the narrow band and the slip wall are refined.}
\label{fig:setup_bosss}
\end{figure}

    \paragraph{Numerical setup:} Since the considered case is symmetric, the computational domain is restricted to the lower side of the capillary, see Figure \ref{fig:setup_bosss}. The boundary condition on the symmetry plane, i.e. the upper boundary in figure \ref{fig:setup_bosss}, is given by a free-slip condition with a contact angle of $\theta_e= 90^\circ$. According to the reference model, which assumes a static circular shape of the interface, we use the static contact model formulation of \eqref{eq:GNBC} with $\fcL = 0$ for the slip wall. At both ends of the capillary we employ an inflow and outflow boundary condition, respectively, with $\mu\left( \partial_n \boldsymbol{u} + \nabla u_n \right) - p \boldsymbol{n}_{\partial \domain} = 0$.
For the initial position of the interface the same start-up procedure is employed as described in \ref{sub:res_ALE} for the ALE method yielding the setup depicted in Figure \ref{fig:setup_bosss}. 
For all the following simulations we use a polynomial degree of $k=2$ for the velocities and $k-1$ for pressure. The grid cells are equidistant and are defined by the grid size $h$. 

\paragraph{Convergence study:} The considered grid sizes for this study are $h = [R, R/2, R/4, R/8, R/16, R/32]$. For the first five sizes unrefined base grids are chosen. For $h = R/16, R/32$ a base grid of $h=R/8$ is used with adaptive mesh refinements on level 1 and 2. On level 1 a base cell is refined in four equidistant cells and into 16 cells on level 2. The cells to be refined are those located in the narrow band around the interface and at the slip wall, see figure \ref{fig:setup_bosss}. Thus the total number of cells ranges from 8 to 2048. Hence, the corresponding number of degrees of freedom (NDOF) for a polynomial degree of $p=2$ is calculated by $15 (K + K_\text{near})$, where $K$ denotes the number of grid cells and $K_\text{near}$ the number of near-band-cells. The additional $15 K_\text{near}$ DOFs result from the separate DOFs of the second phase. For the convergence study the time step size is set constant for all runs and is chosen based on the smallest grid size. 
\begin{figure}[H]
	\centering
	\subfigure{\includegraphics[width=0.48\textwidth]{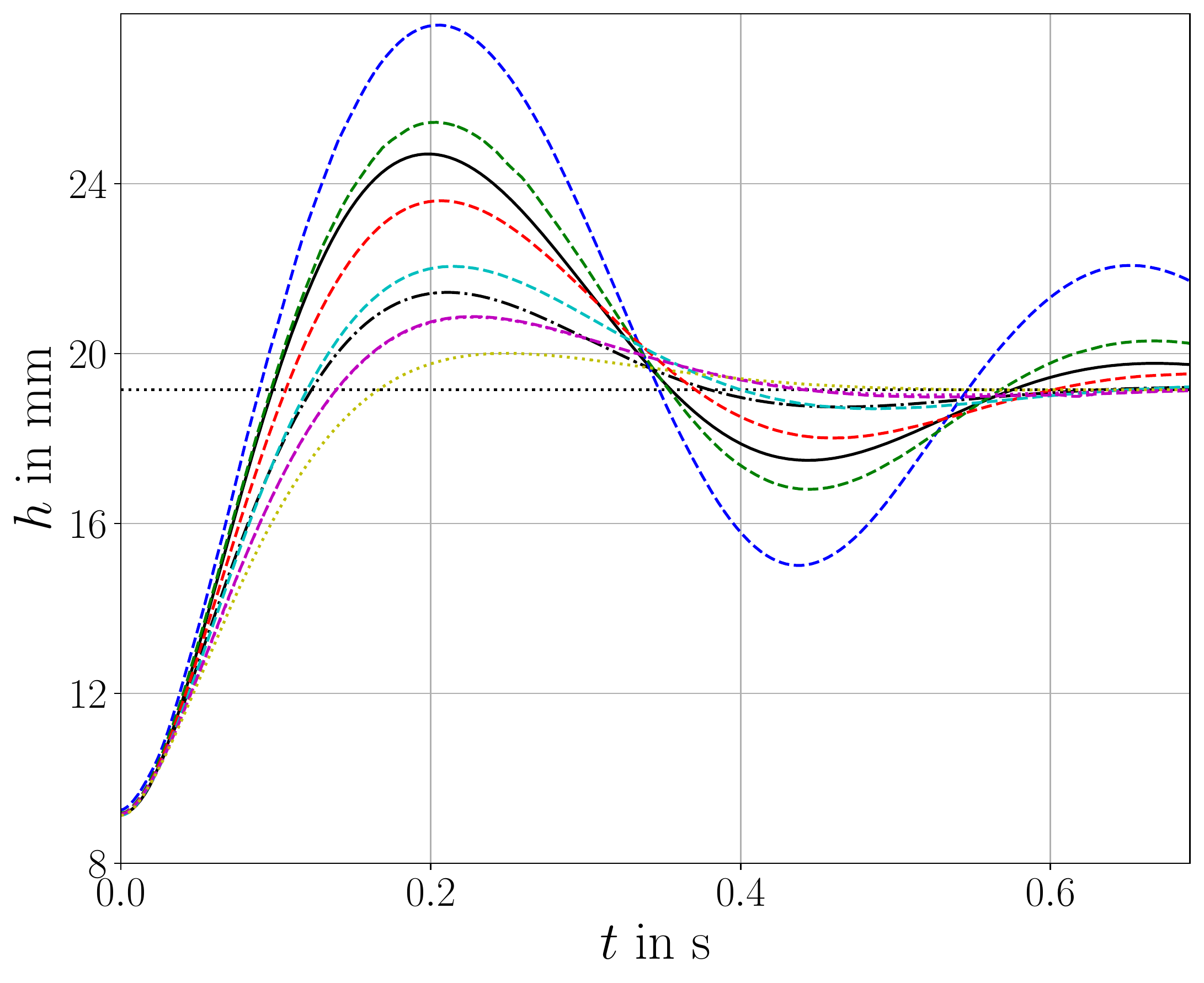}}
	\subfigure{\includegraphics[width=0.48\textwidth]{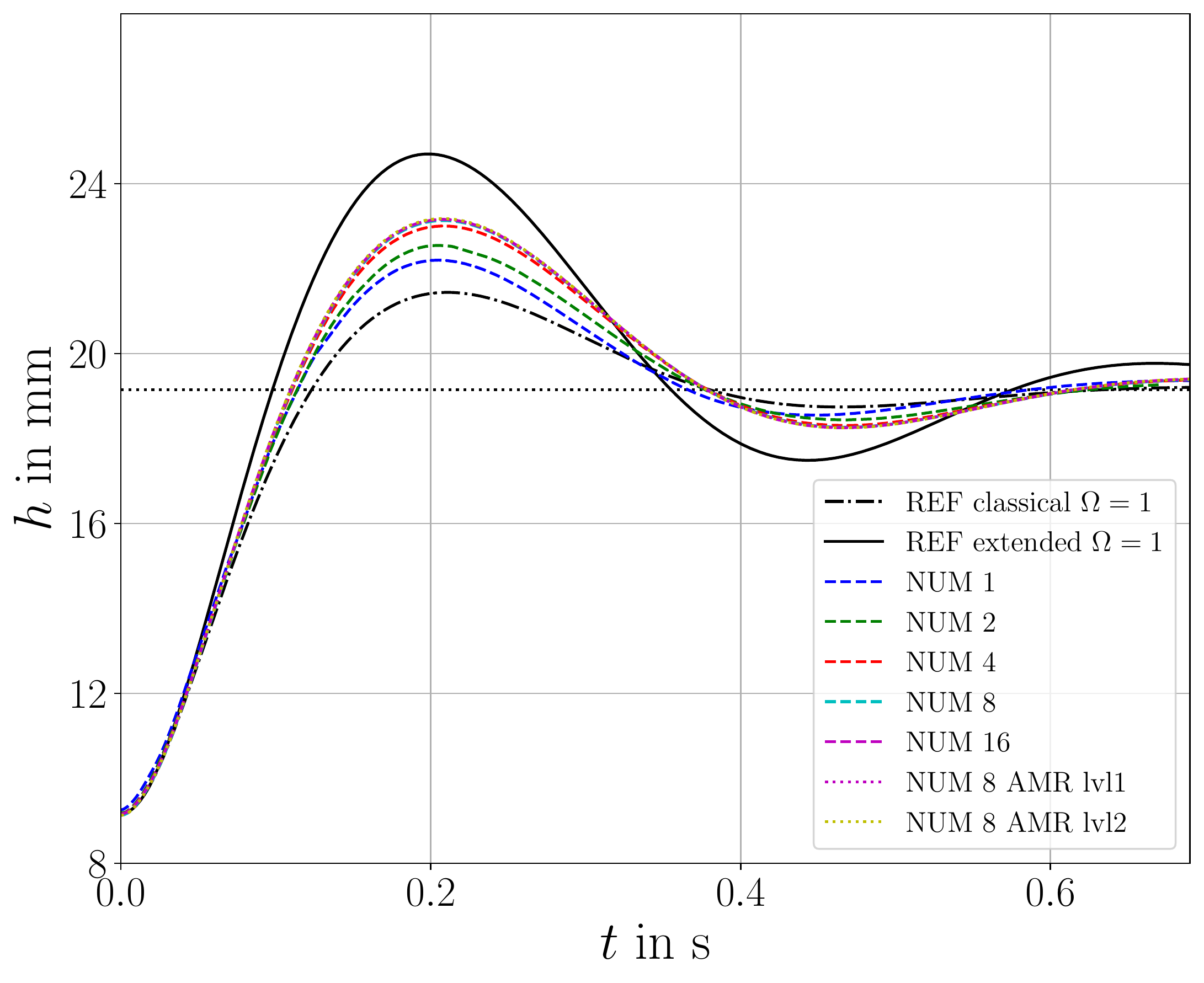}}
	\caption{Mesh study with BoSSS using  an effective slip in order of the grid size (left) and a resolved Navier slip boundary condition (right).}
	\label{fig:mesh_study_bosss}
\end{figure}

Since the discretization of the GNBC \eqref{eq:GNBC} in XDG does not need to manipulate/prescribe the contact line velocity and contact angle directly, it does not allow a numerical slip as for the other methods. But, as described above, the numerical slip introduces an effective slip length, which is in order of the given grid size. In order to compare the effect of a numerical slip approach for XDG, the prescribed slip length in the mesh study is set to the corresponding grid size. The results are shown on the left hand side in figure \ref{fig:mesh_study_bosss}. As observed before, no converged solution is reached. Furthermore, the solution vary more and show much higher oscillations for the coarser meshes.    
 On the right hand side of Figure \ref{fig:mesh_study_bosss},
a resolved Navier slip length of $L = R/5$ is used. As with the other approaches, the solution is mesh converged. Compared to the reference solutions, the results lie between the classical no-slip and the extended slip solution. However, the numerical results on the coarser meshes converge in contrast to the other numerical methods from below. This may indicate, a stronger influence of the slip/contact line dissipation due to the polynomial solution approximation in the cells.

    \subsection{Comparison}\label{sub:comparison}
In this Section we consider on one hand a code to code comparison between the four numerical approaches and on the other hand a comparison to the aforementioned 1D reference models. Firstly, the results from the previous convergence studies for $\Omega = 1$ are compared between all numerical implementations. Secondly, results for all implementations for a reduced slip length are discussed to illustrate the influence of a reduced slip length. At last, the extended reference solution and the results of the $\Omega$-study are compared using a scaling that is available in literature.

\subsubsection{Convergence study and varying slip length}

\begin{figure}[H]
  \centering \subfigure{\includegraphics[width=0.48\textwidth]{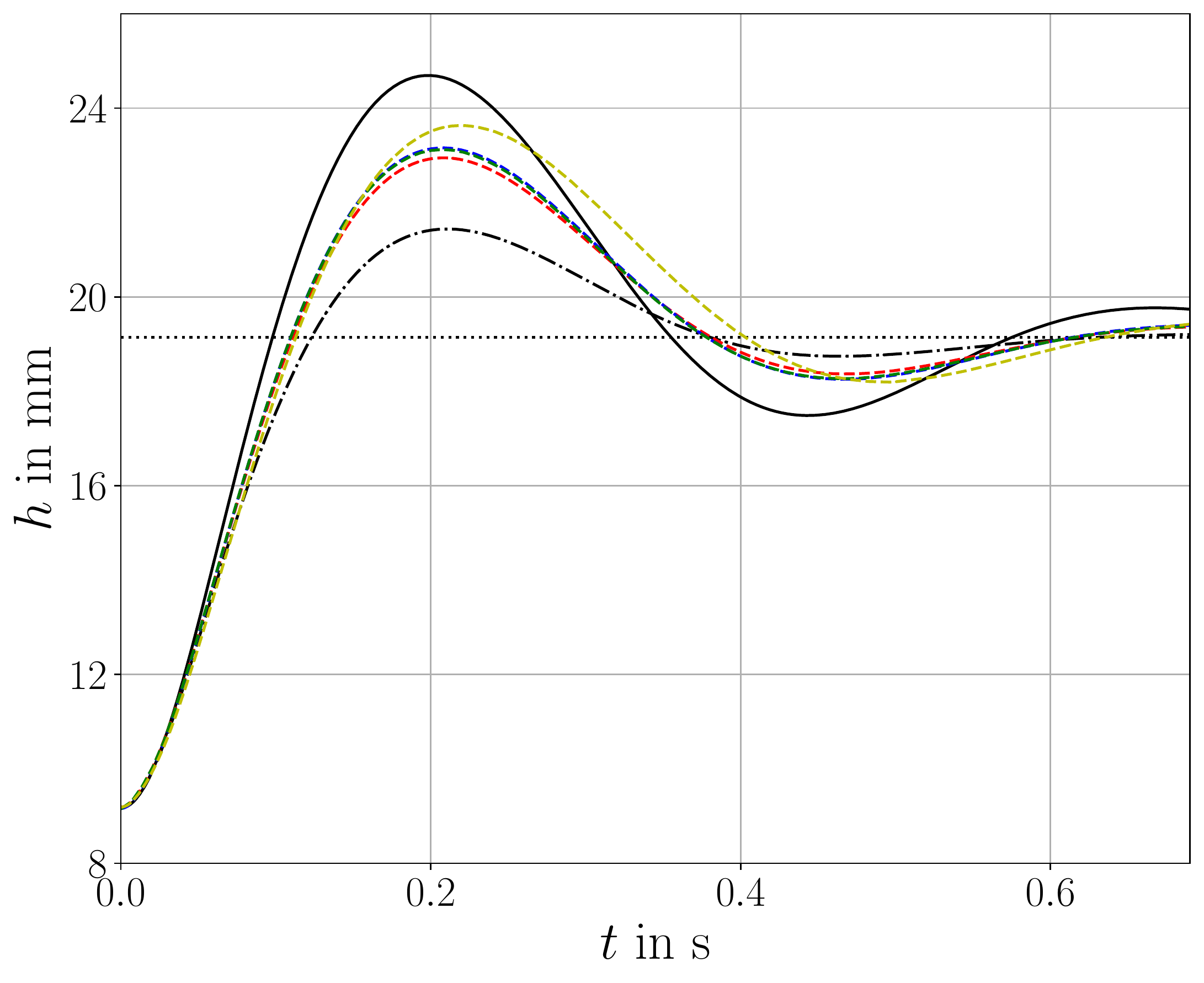}}
	\subfigure{\includegraphics[width=0.48\textwidth]{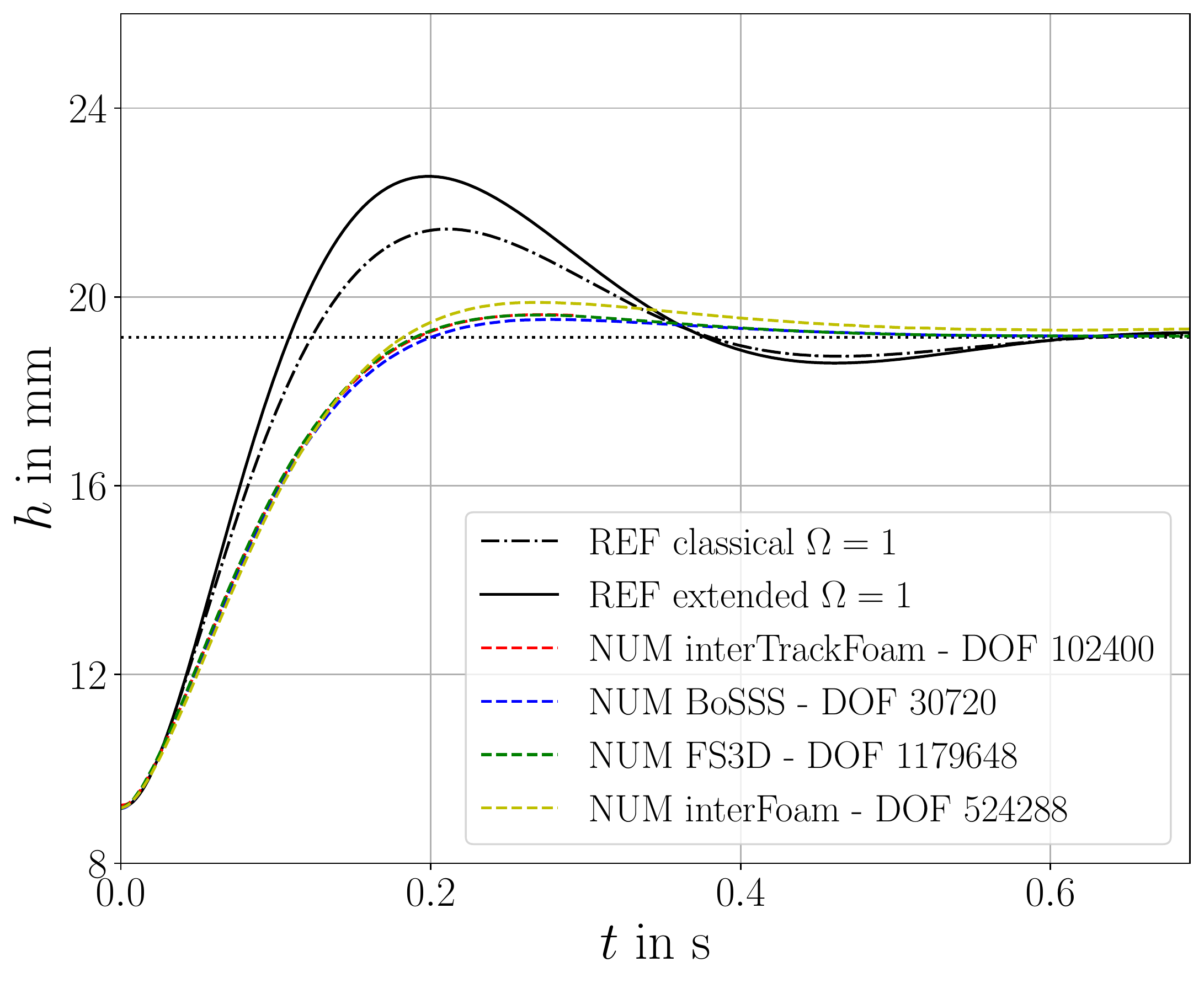}}
	\caption{Comparison of the convergence study for $\Omega = 1$ with different slip lengths: (left) slip length $L = R/5 = \SI{1}{\milli \meter}$ and (right) slip length $L = R/50 = \SI{0.1}{\milli \meter}$.}
	\label{fig:conv_study}
\end{figure}

In the previous sections two convergence studies for each implementation have been shown. All methods show a mesh dependence for the no-slip/numerical slip approach. For a resolved slip length of $L = R/5 = \SI{1}{\milli \meter}$ all methods, except for the algebraic VOF method (interFoam), yield mesh converged results. Figure \ref{fig:conv_study} shows a comparison of the finest solution from each method for the case with Navier slip boundary conditions and a selected interFoam solution that fits best the results of the other three solutions. The overall agreement is good. BoSSS and FS3D results are in excellent agreement. The oscillation amplitude of the ALE method is slightly smaller. In comparison, the solution computed with interFoam oscillates with a slightly stronger amplitude as well as a shorter frequency. 
Comparing the numerical solutions with the classical no slip solution \eqref{eq:ode_classical} and the extended version \eqref{eq:ode_extended}, it can be seen that all solutions exhibit qualitatively the same dynamic behaviour in this case with a larger slip length. The classical no-slip solution shows a less and the extended solution a more dynamic behavior then the full continuum solutions. As this seems to be a trend for all considered values of $\Omega$ with $L=R/5$, (see next section) we consider the case of a reduced slip length with $L = R/50 = \SI{0.1}{\milli \meter}$ which is shown on the right of Figure \ref{fig:conv_study}. With this reduced slip length, the full continuum mechanical results are less dynamic than both reference solutions. The direct cause is unclear as the reference solution make various approximations. For example, their derivation assumes a quasi-stationary Poiseuille flow regime over the most part of the capillary and fully neglects all dissipation effects caused by the velocity field in contact line vicinity. However, most surprising is the fact that the full numerical solution does barely show any oscillations even though the slip length is still far away from a nano-scale slip length reported in literature. Hence, reducing the slip length also reduces the dynamic of the rise similar to the behavior for the numerical slip approach. Yet, using Navier slip boundary condition still gives convergent results. This observation has already been made in \cite{Gruending2019} for the ALE approach.

\subsubsection{$\Omega$-study}
\label{sec:NumResults_Ostudy}
\label{subsub:OmegaStudy}

\begin{figure}[htb]
	\centering
	\subfigure{\includegraphics[width=0.48\textwidth]{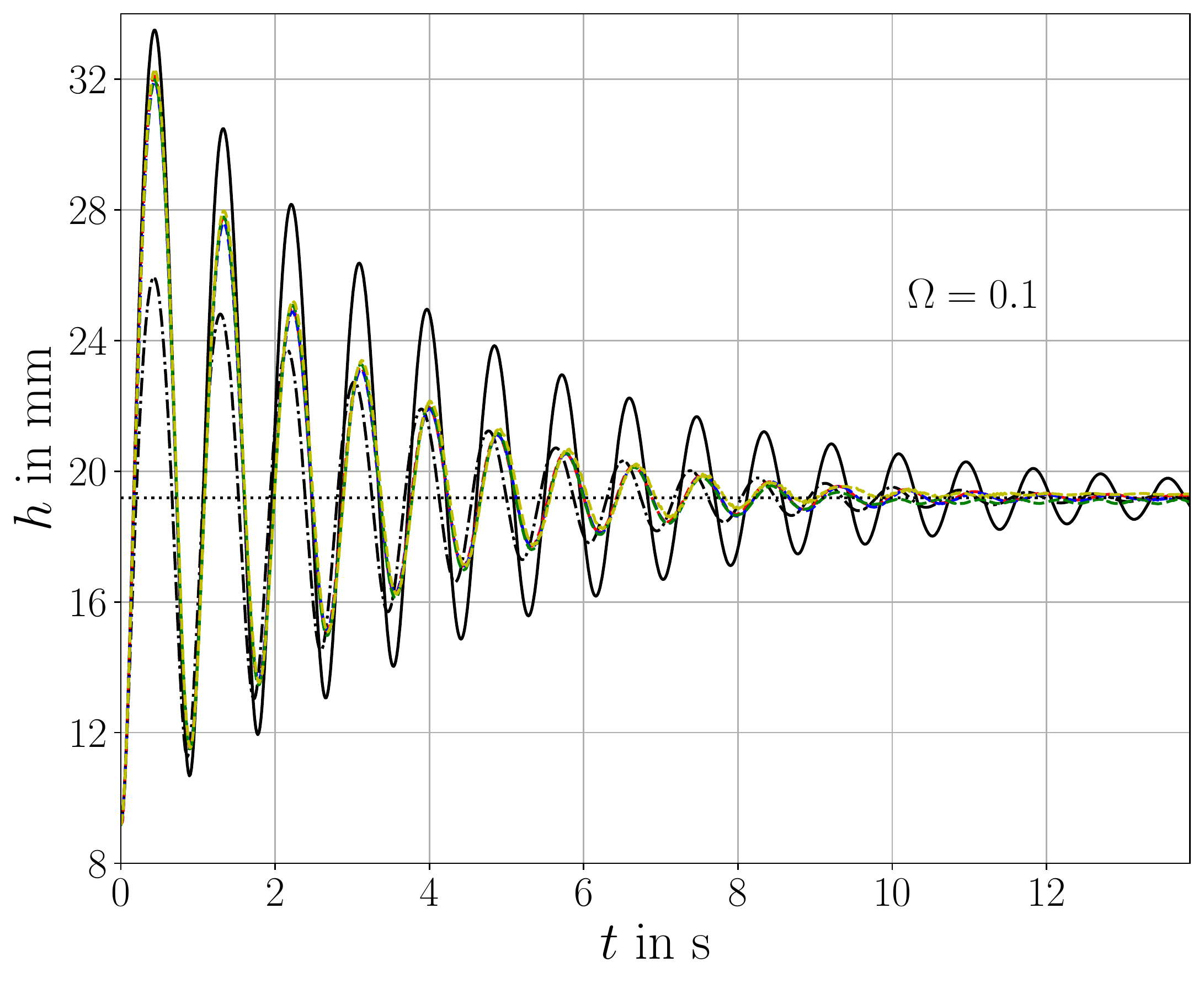}}
	\subfigure{\includegraphics[width=0.48\textwidth]{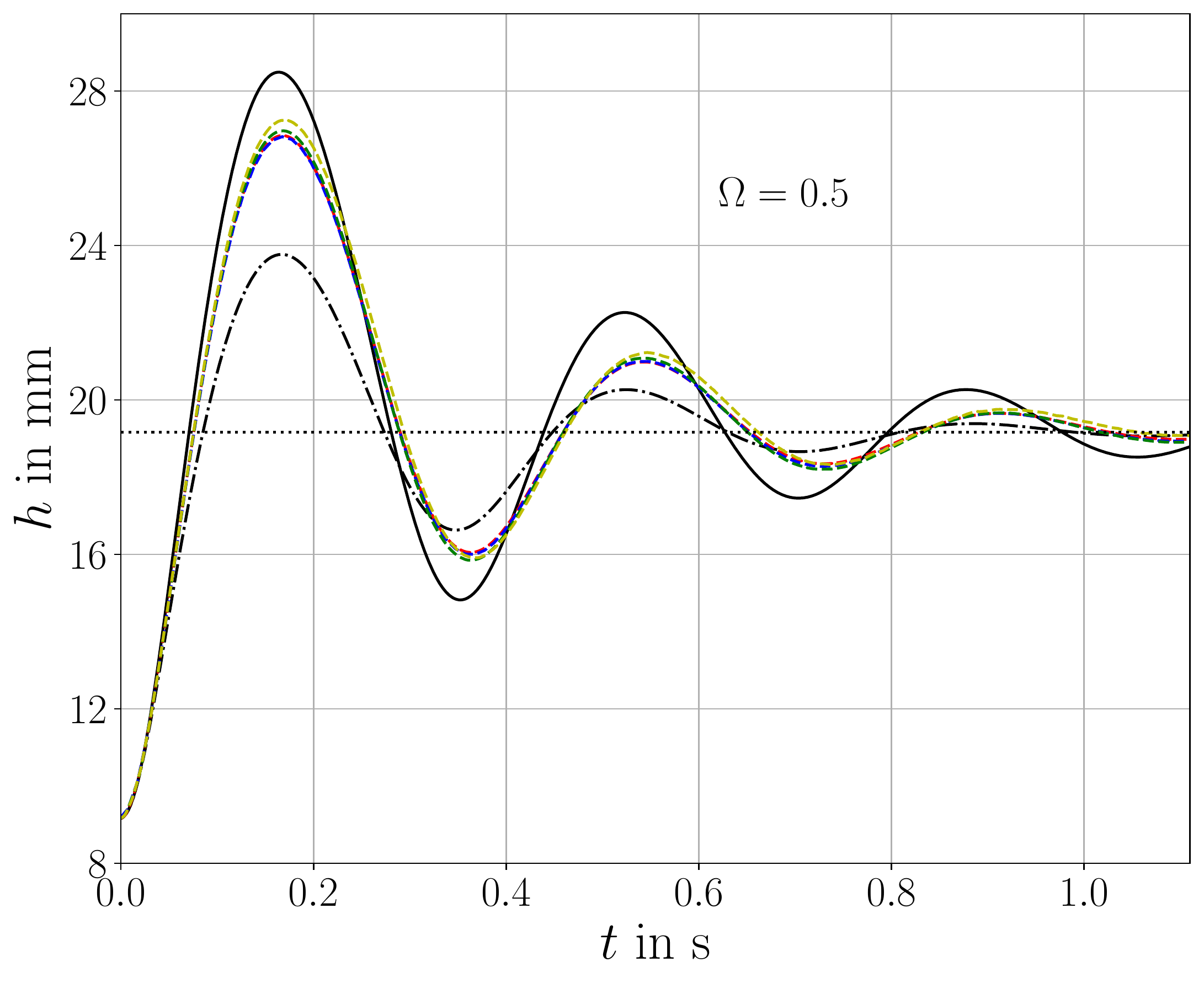}}
	\subfigure{\includegraphics[width=0.48\textwidth]{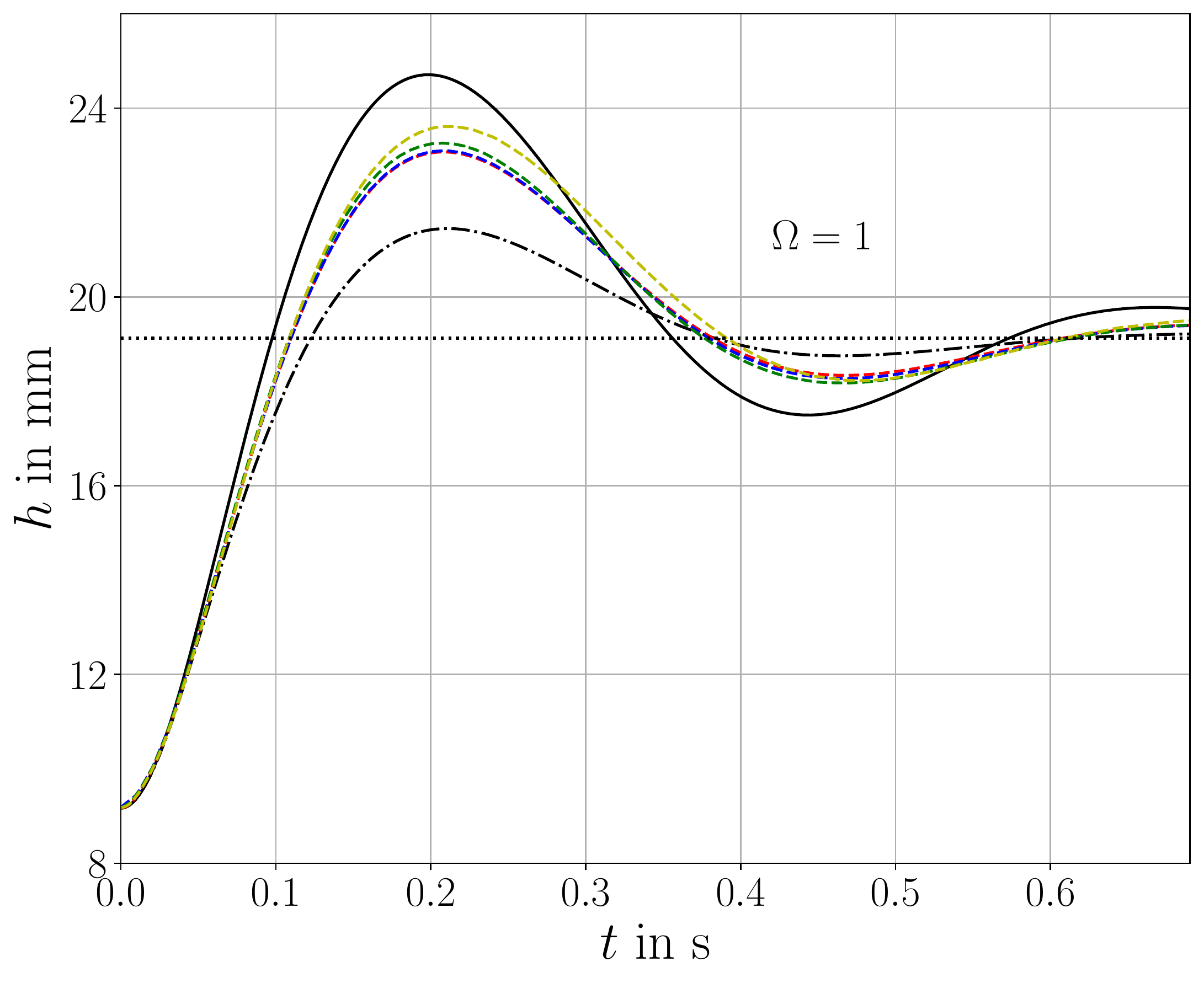}}
	\subfigure{\includegraphics[width=0.48\textwidth]{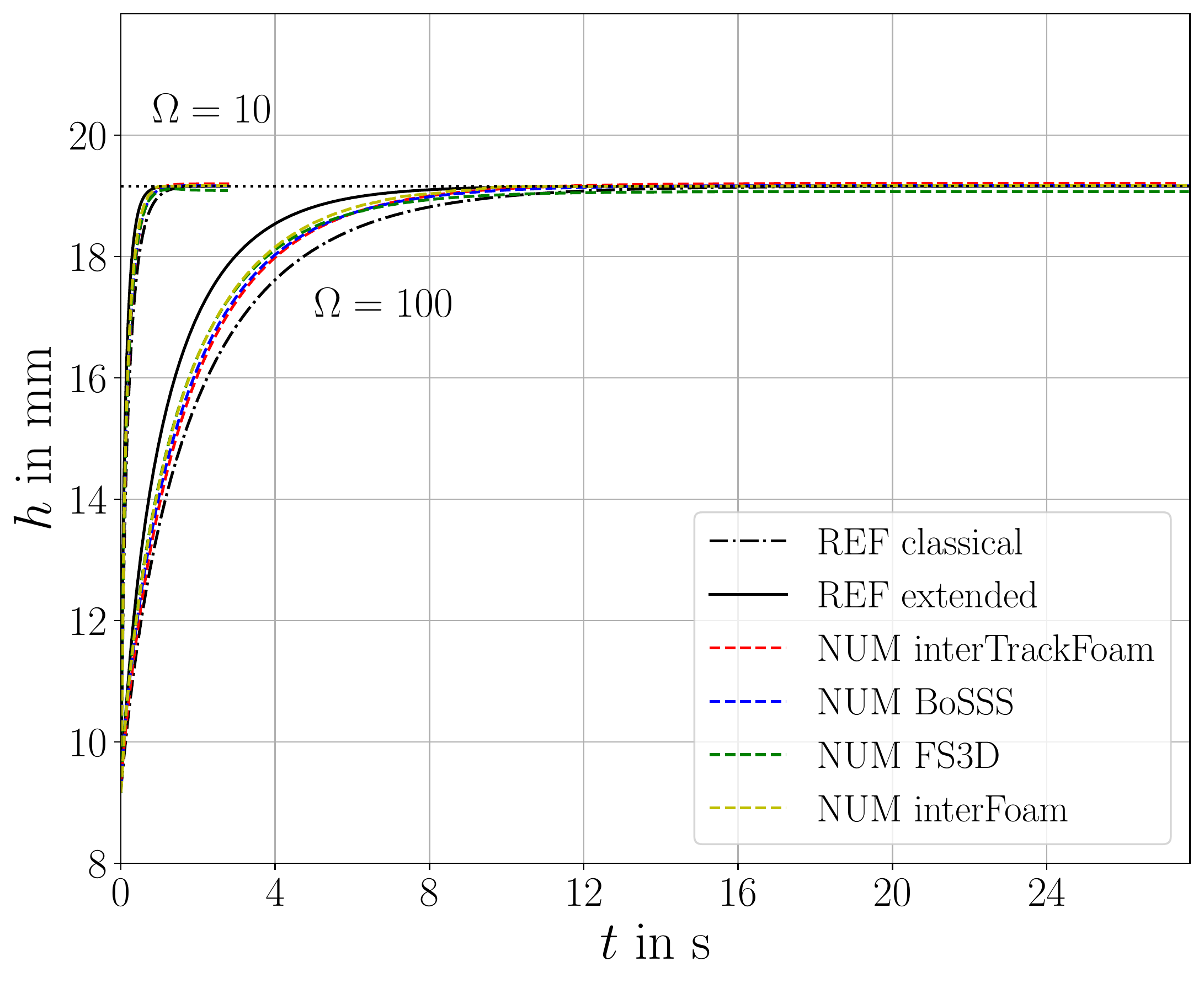}}
	\caption{Unscaled solutions of all numerical methods and reference models for the $\Omega$-study: $\Omega = 0.1$ (top left), $\Omega = 0.5$ (top right), $\Omega = 1$ (bottom left) , $\Omega = 10, 100$ (bottom right).}
	\label{fig:omega_study_unscaled}
\end{figure}

\paragraph{Unscaled results}
Figure \ref{fig:omega_study_unscaled} shows the results for the capillary rise problem obtained with the four different numerical methods using a slip length of $L=R/5$ and a comparison with the results for the 1D ode models. The results obtained with the different numerical methods show good to excellent quantitative agreement compared to each other throughout all simulations. All solutions show a strong oscillation that is damped over time and levels at the stationary height predicted by \eqref{eq:stationary_rise_height}.

For $\Omega=0.1$, shown in the top left graph, all numerical solutions coincide up to approximately $t=\SI{5}{\second}$. From then on, first the interFoam solution, followed by the results obtained with FS3D start to show minor deviations from the remaining ALE and BoSSS results. In addition, the period of oscillation shows a good agreement during the initial oscillations up to approximately $t=\SI{5}{\second}$ and the phase shift increases over time. 

The first maximum is well met by the extended model. However, this reference model shows far less dampening than the full continuum solution. 
The first maximum is underpredicted by the classical model by about \SI{50}{\percent}. After roughly four oscillation periods the continuum and classical model show a good agreement and a similar dampening behavior up to approximately \SI{9}{\second} where agreement between the different continuum solutions is lost. Again, the continuum and the reference  solution are out of phase at this point. With increasing $\Omega$ the solutions show less oscillations as can be seen from the results for $\Omega=0.5$ to $\Omega=100$ with nearly aperiodic behaviour for $\Omega=1$ as shown in the bottom left graph of figure \ref{fig:omega_study_unscaled}. For $\Omega = 0.5$ and $\Omega=1$ interFoam and (less pronounced) FS3D over and undershoot the maximum and minimum of the ALE and BoSSS solutions. 

A general trend for this slip length is that the continuum solutions show a dynamic that is somewhat between the two reference solutions for all values of $\Omega$. Furthermore, all approaches give the same stationary height in agreement with the expected value according to equation \eqref{eq:stationary_rise_height}. The cases for $\Omega=10,100$ do not show any oscillations. Overall, the solution behavior shows a stronger dependence on $\Omega$ for decreasing values of $\Omega$.

\begin{figure}[H]
	\centering
	\subfigure{\includegraphics[width=0.48\textwidth]{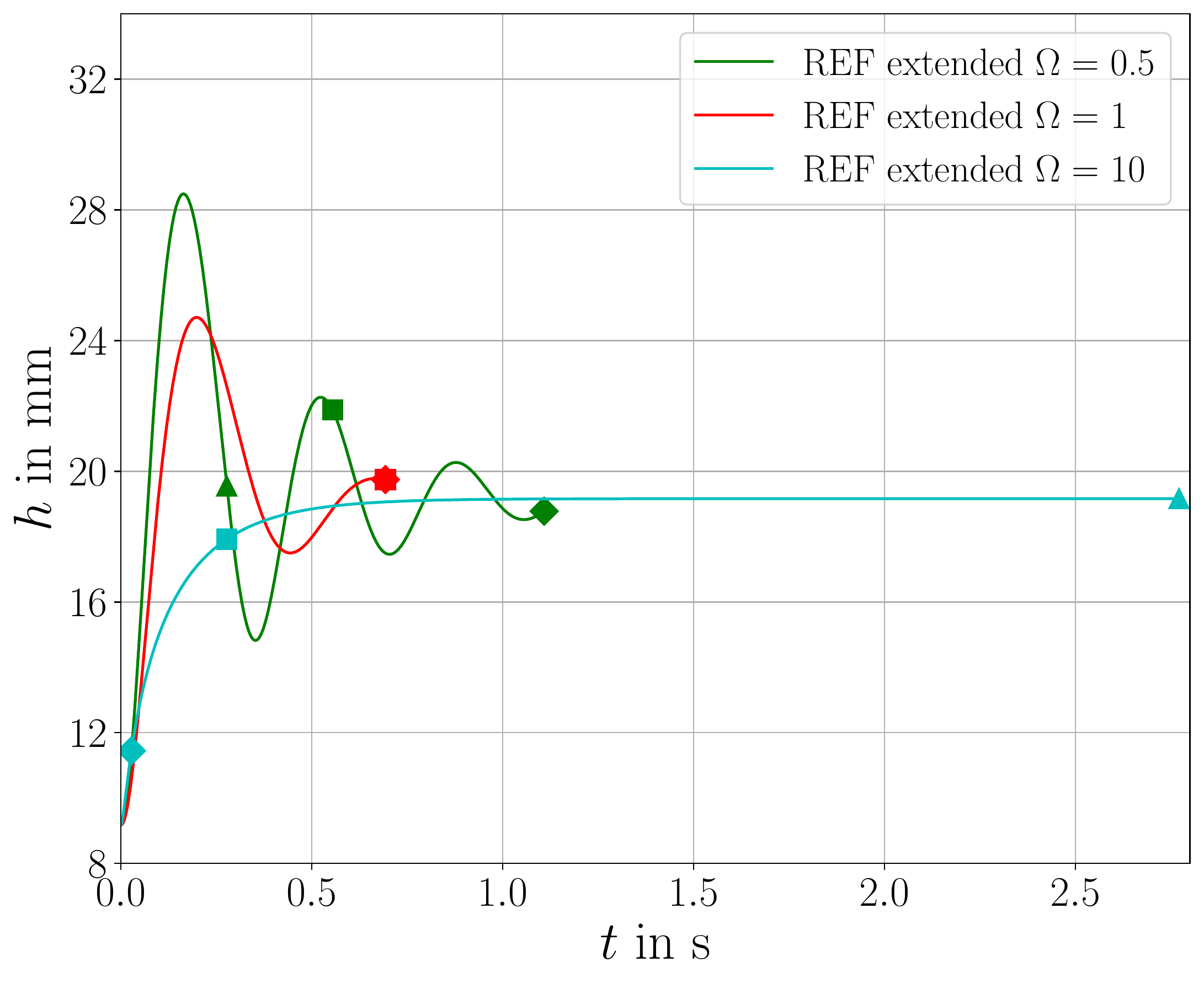}}
	\subfigure{\includegraphics[width=0.48\textwidth]{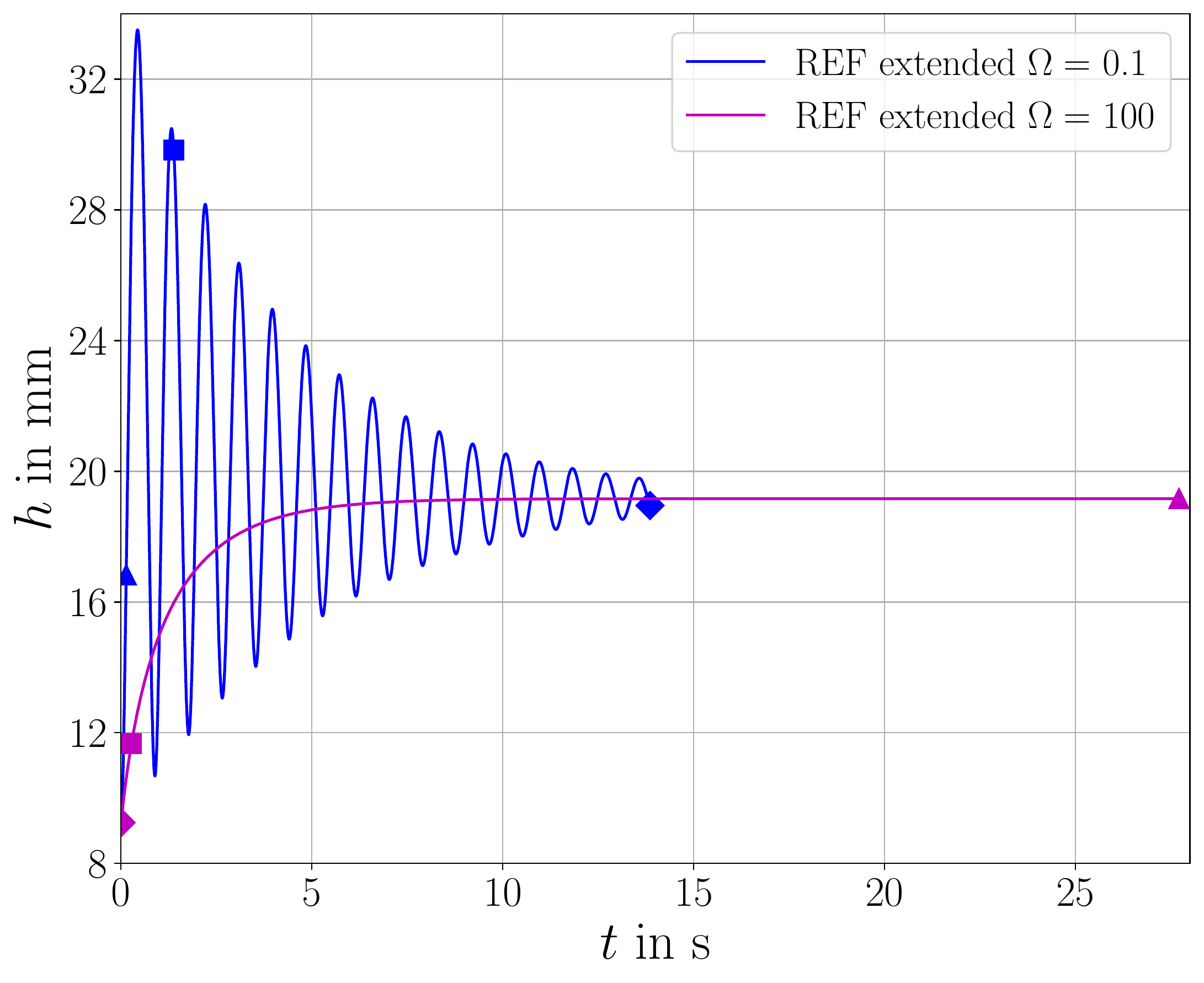}}
	\subfigure{\includegraphics[width=0.48\textwidth]{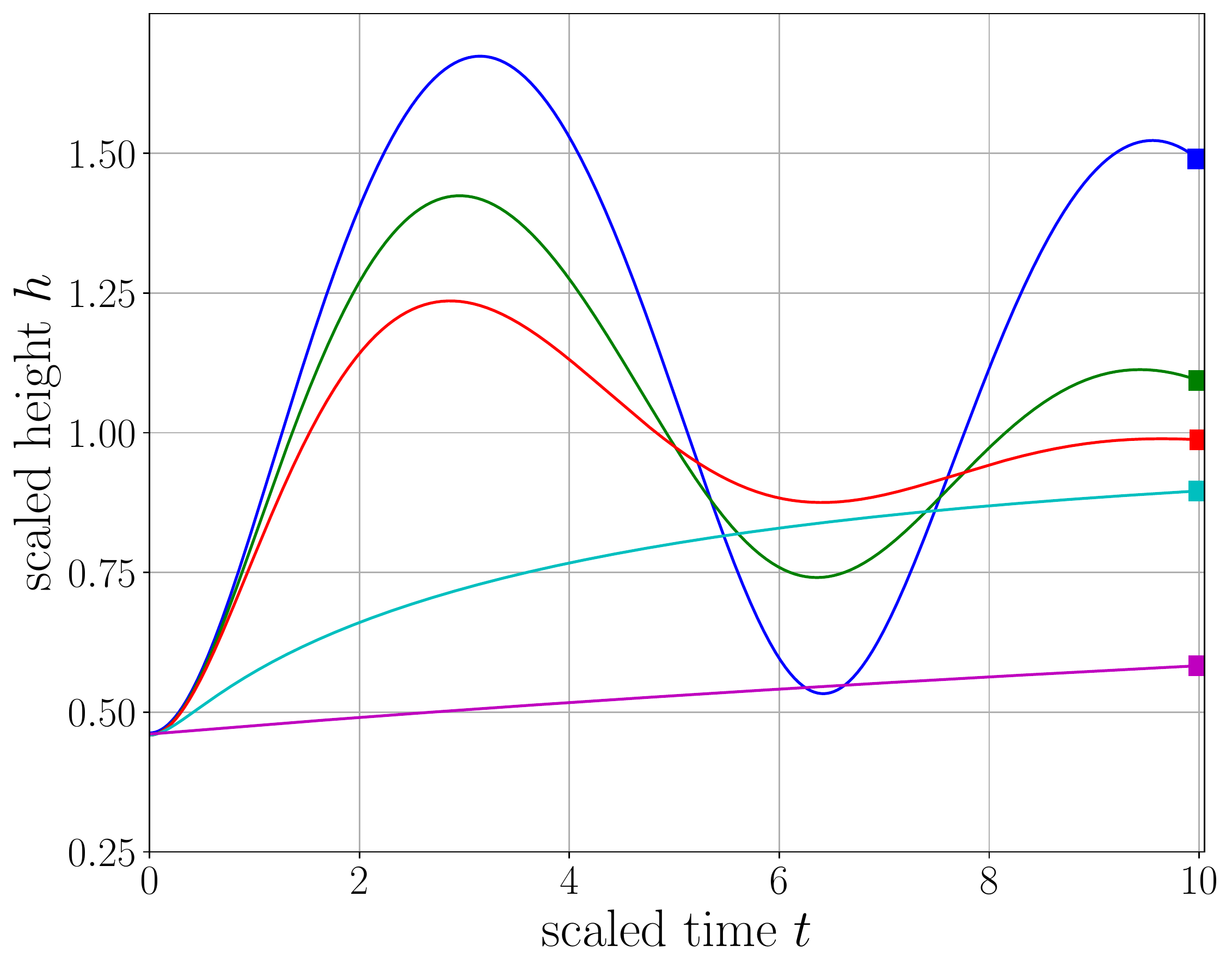}}
	\caption{First row: unscaled extended reference solutions for the $\Omega$-study.  Second row: Scaled solutions of the extended reference model for the $\Omega$-study using scaling II. The symbols mark the end times for each individual scaling: triangles for scaling I, squares for scaling II and diamonds for scaling III.}
	\label{fig:omega_study_ref_scaled}
\end{figure}

\paragraph{Scaled results}
To provide generalized benchmark data, the results for the $\Omega$-study are scaled according to the discussion given in section \ref{sub:validation_case}. Here, we also consider scaling II as it allows to visualize the results of $\Omega=10$ and $\Omega=100$ without overlap. This scaling uses inertia and gravity as scaling forces with viscosity as basic parameter. This means that an increasing value of $\Omega$ corresponds to an increasing influence of viscosity, cf. \cite{Fries2009}. As the approach and results are similar for the other scalings, we will only consider scaling II here. For completeness, the remaining scalings I and III can be found in the Appendix~\ref{app:comparison}.
An example how the scaling influences the visual representation can be seen in Figure \ref{fig:omega_study_ref_scaled}. 
The top row shows the unscaled results together in two plots; for $\Omega = [0.1, 1, 100]$ in the left and $\Omega = [0.5,10]$ in the right plot. The solutions are plotted until $t_\text{end}$, given in table \ref{tab:physicalParams}, which is the time to cover all scaling given in \cite{Fries2009}. Notice the square symbols on each graph. Applying scaling II gives the lower plot in Figure \ref{fig:omega_study_ref_scaled} and ensures that these squares are all located at the scaled time $t=10$. The same holds for scalings I and III which are depicted by triangles and diamonds, respectively.

\begin{figure}[H]
	\centering
	\subfigure{\includegraphics[width=0.48\textwidth]{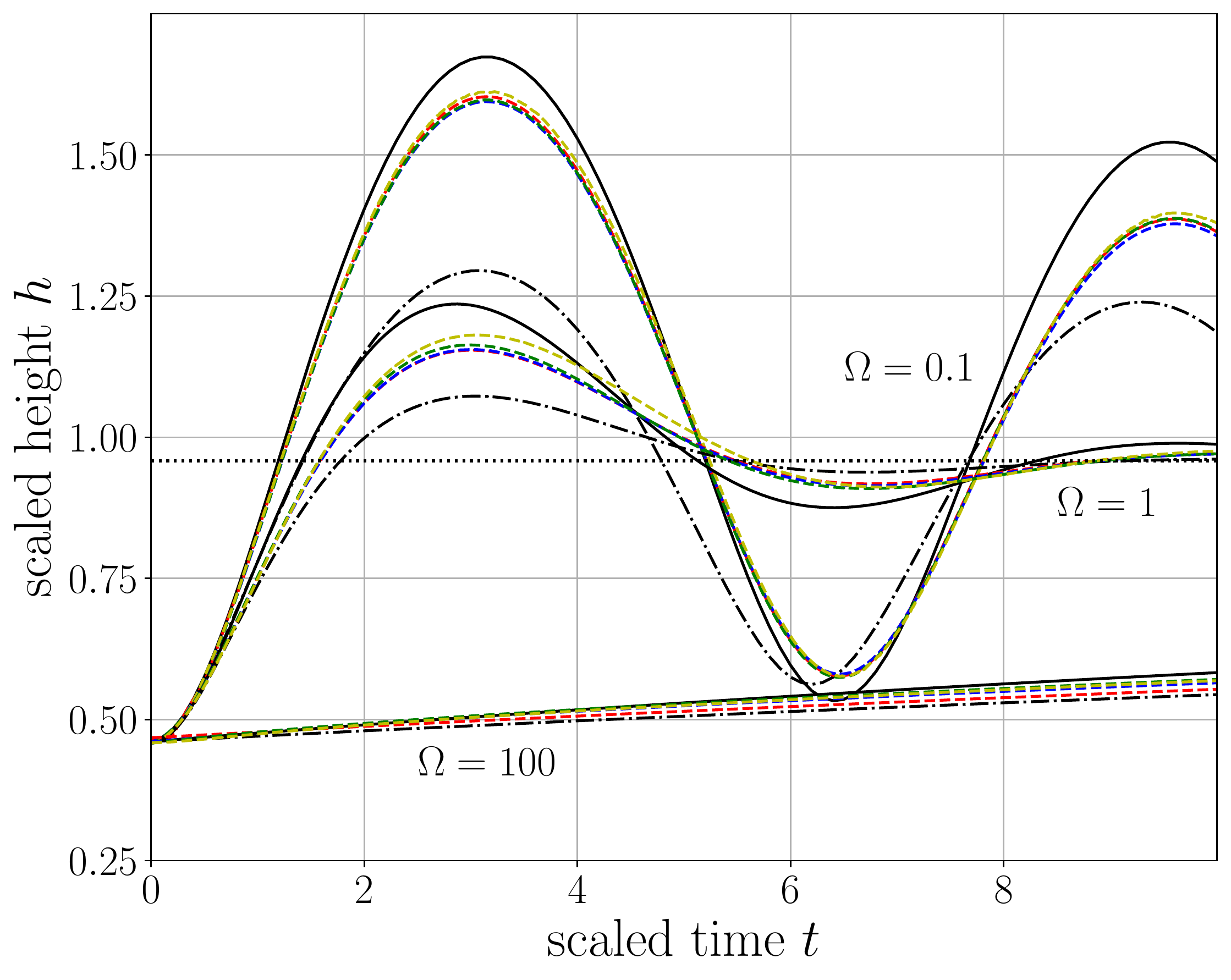}}
	\subfigure{\includegraphics[width=0.48\textwidth]{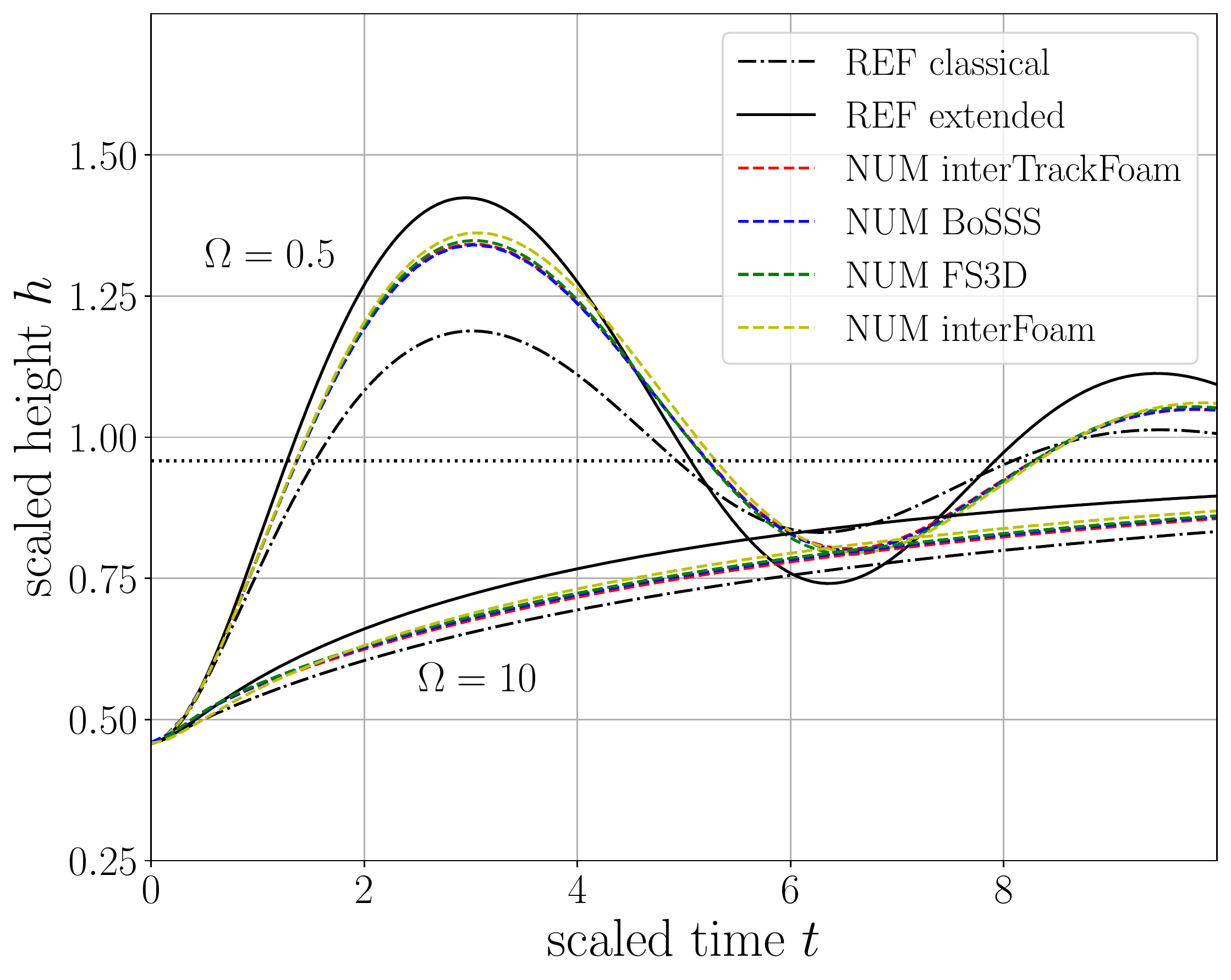}}
	\caption{Dynamic behavior for the $\Omega$-study with scaling II. left: $\Omega = [0.1, 1, 100]$, right: $\Omega = [0.1, 10]$.}
	\label{fig:omega_study_compar_Scale2}
\end{figure}

Applying scaling II to the the solution of the continuum mechanical problem from all four numerical approaches together with the two reference solutions are shown in Figure~\ref{fig:omega_study_compar_Scale2}. The results for $\Omega = [0.1, 1, 100]$ are displayed in the left and $\Omega = [0.5,10]$ in the right plot. 

As for the unscaled results before, all numerical solutions show good to excellent agreement among each other with slight deviations for interFoam in the cases $\Omega = 0.5$ and $\Omega = 1$. For $t < 2$ the results from the four numerical methods show an excellent agreement. Yet, the scaling may slightly increase the differences between these results and the deviations are negligible in comparison to the deviation from the two reference solutions. 

Again, comparing the numerical solutions to the 1D reference model and the classical solution, the basic dynamical behavior depending on the different $\Omega$ is qualitatively met well. In addition, the qualitative agreement is similar to the unscaled results. The solutions for $\Omega = 0.1$ and $\Omega = 0.5$ show the largest deviations, whereas for larger $\Omega$ the agreement is improving. For the small $\Omega = [1, 0.5, 0.1]$ and thus more dynamic cases, the numerical solutions agree well with the point in time of the first peak. However, for longer times, oscillation amplitude and period start to deviate.

In conclusion, the qualitative behavior is similar to the unscaled results. Furthermore, there does not seem to be a pronounced viscous, gravitational or inertial influence (see Appendix \ref{app:comparison}) that can be obtained from the scaled results. This concerns both, the agreement for the code to code as well as the model to code comparison.

\section{Conclusion}
Comprehensive numerical simulations for the rise of liquid between two plates has been performed comparing the results between four numerical methods and two ODE models. The methods include an ALE approach (interTrackFoam), an algebraic (interFoam) as well as a geometric (FS3D) Volume of Fluid approach, and a level-set based extended discontinuous Galerkin discretization (BoSSS). The results have been compared to the classical model for liquid rising in a capillary and an extended reference solution that includes the influence of a Navier slip boundary condition at the capillary walls. The available theory for the ODE models has been used to estimate useful parameters for the comparisons and give a non-dimensionalized representation of the results.

The stationary rise height of all four implementations agrees well with the corrected analytic reference for the considered regime of the Eotvos number. Furthermore, for the considered mesh resolutions, a no slip/numerical slip combination does not yield mesh convergent results which, on the other hand, could be obtained with a Navier slip boundary condition. This result concerns all four approaches and illustrates that using numerical slip does not seem to be viable approach in general. Regarding the Navier slip boundary condition, a monotonic mesh convergent solution could be obtained for the ALE, BoSSS and FS3D implementations. The converged solutions of these three codes show excellent agreement for all conducted simulations. Improvements for the interFoam approach were necessary in order to use a Navier slip boundary condition at the wall. Without these adaptations, a smearing of the phase indicator field lead to unusable results. While selected solutions of the considered interFoam solver can achieve a good quantitative agreement with the other three implementations for all considered cases, no mesh converged results could be obtained for the considered cases.

While the extended reference model reduces to the classical model without slip for a vanishing slip length, this is not the case for the continuum solution. Decreasing the slip length in the full continuum mechanical simulations significantly reduces the dynamics of the problem. This does not only change the results quantitatively but also qualitatively as oscillations expected due to the reference model vanish for the continuum solutions with decreasing slip length. 

Using a slip length on the scale of the capillary, the qualitative agreement between the full continuum solutions and the ODE models is good including the onset of rise height oscillations. However, the quantitative agreement is reasonable between both reference models and the full solution of the continuum mechanical problem and decreases for decreasing $\Omega$, i.e., we observed an increasing agreement for cases with weak and absent oscillations. This behavior can most likely be attributed to the approximate nature of the reference models. To provide comparable benchmark data, the obtained results are suitably scaled and made available online. 

\section{Acknowledgement}
We kindly acknowledge the financial support by the German Research Foundation (DFG) within the Collaborative Research Centre 1194 ``Interaction of Transport and Wetting Processes'', Project B01, B02, and B06. Calculations for this research were conducted on the Lichtenberg high performance computer of the TU Darmstadt.

\newpage

\appendix 
\section{Appendix}
\label{app:comparison}
\subsection{Scaling I - the influence of inertia}

\begin{figure}[H]
	\centering
	\subfigure{\includegraphics[width=7cm]{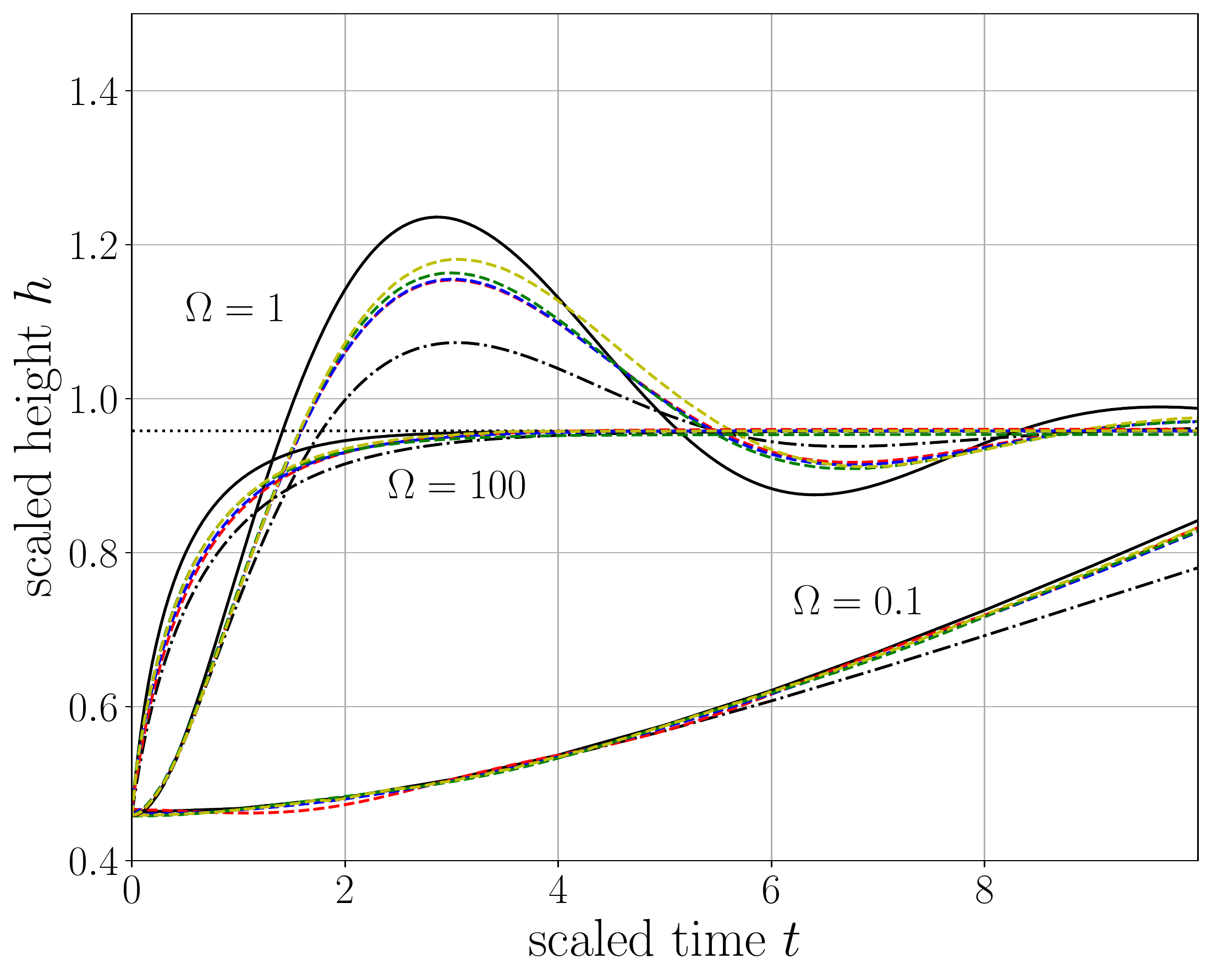}}
	\subfigure{\includegraphics[width=7cm]{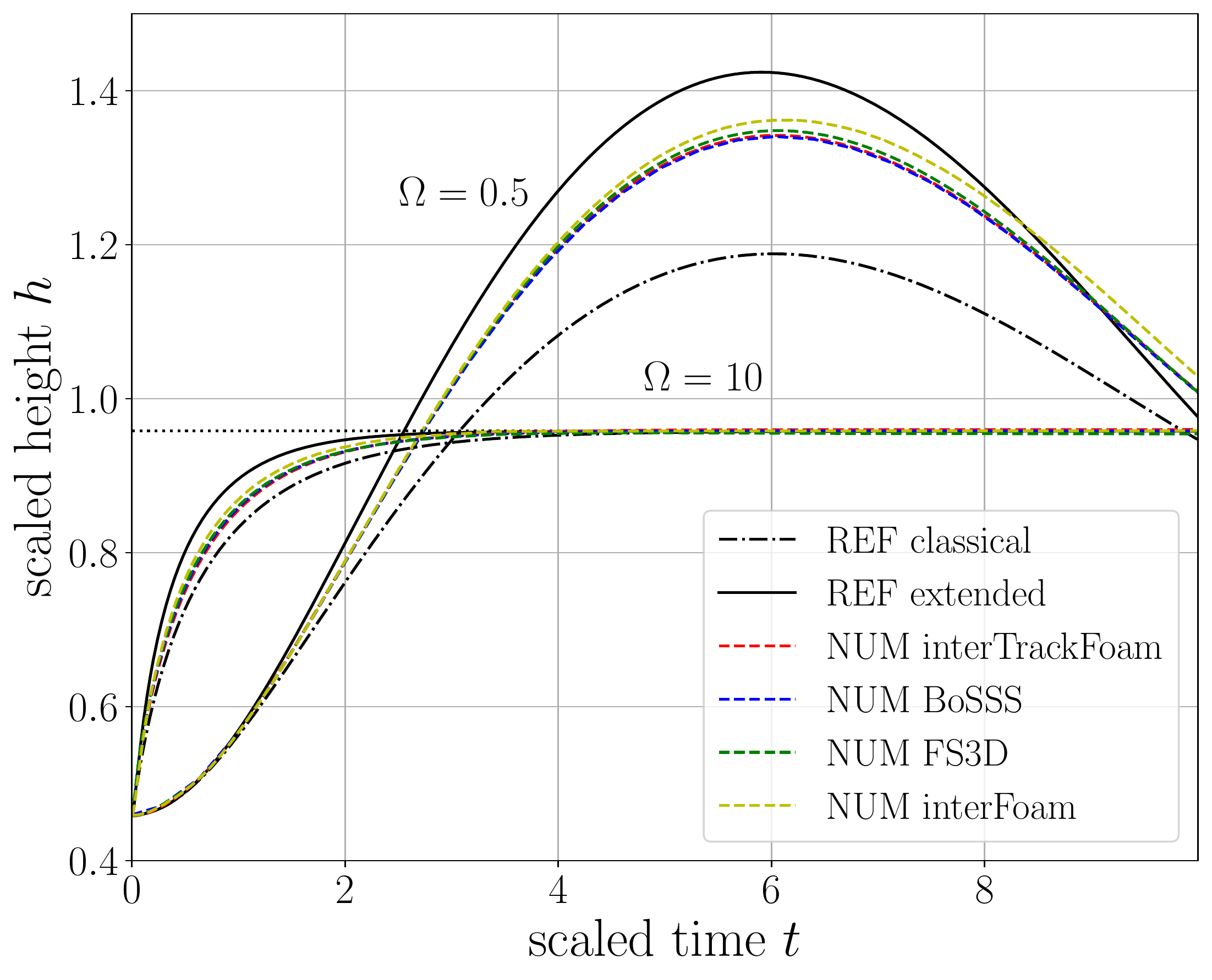}}
	\caption{Dynamical behaviour for the $\Omega$-study with scaling I: left: $\Omega = [0.1, 1, 100]$, right: $\Omega = [0.1, 10]$}
	\label{fig:omega_study_compar_Scale1}
\end{figure}

Figure \ref{fig:omega_study_compar_Scale1} shows the results of the continuum mechanical model for the different numerical methods in comparison with solutions of the 1D reference model for scaling I. With this scaling, an increasing value of $\Omega$ corresponds to an increasing influence of inertia. This influence can, for example, be observed with the results for $\Omega=0.1$ where the results are scaled such that only the initial time before any observation occurs is shown. The general qualitative behavior is the same as for the unscaled results.

\subsection{Scaling III - the influence of gravity}
The third scaling III emphasizes the influence of gravitational forces measured by $\Omega$. Thus, it can be observed, that with decreasing $\Omega$ the amplitudes of the oscillations increase. Furthermore, the scaling gives varying initial heights. Note that this effect is different from the results in \cite{Fries2009} as we used a non-zero initial height where a comparison is given for the classical model with zero as initial height.

\begin{figure}[H]
	\centering
	\subfigure{\includegraphics[width=7cm]{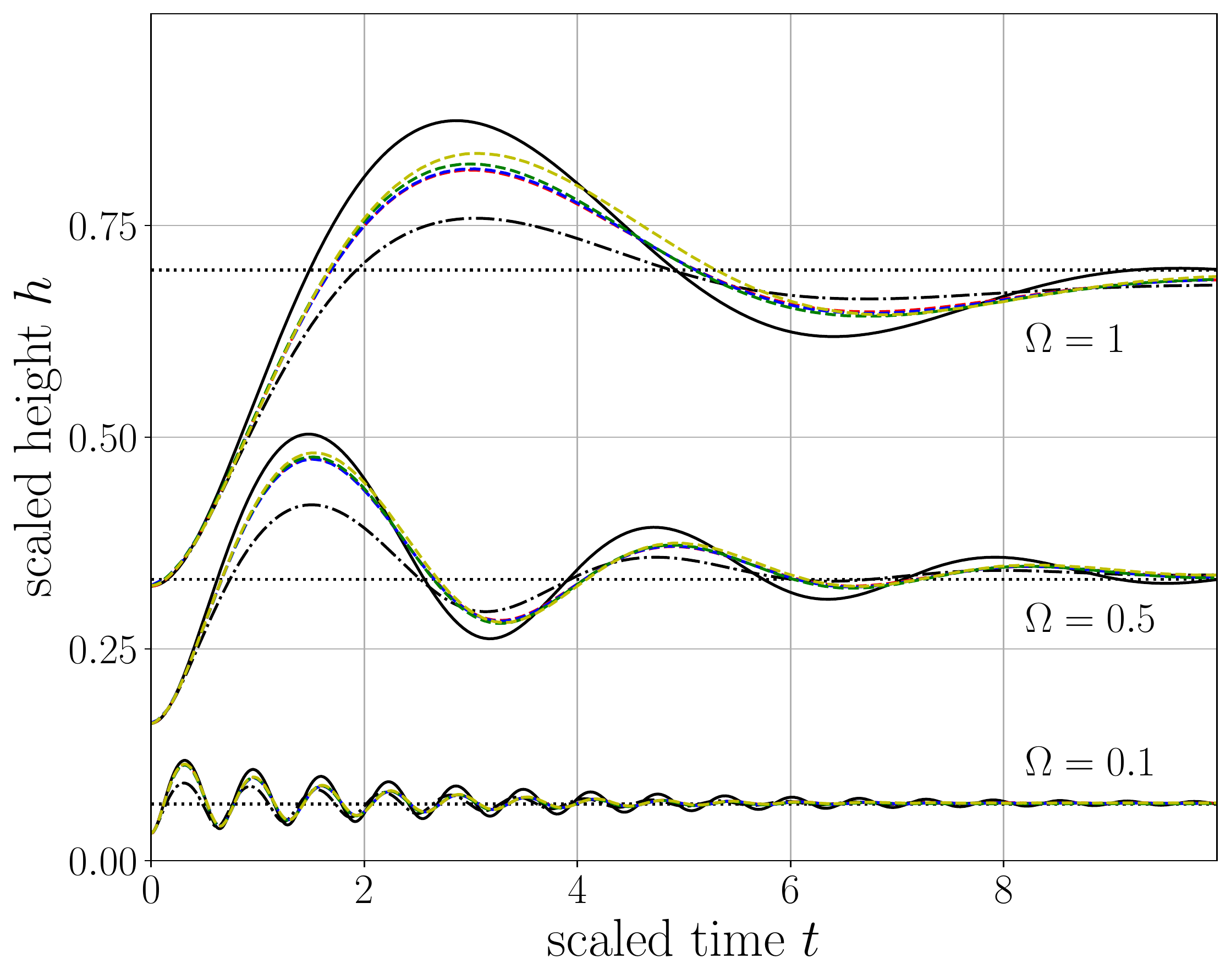}}
	\subfigure{\includegraphics[width=7cm]{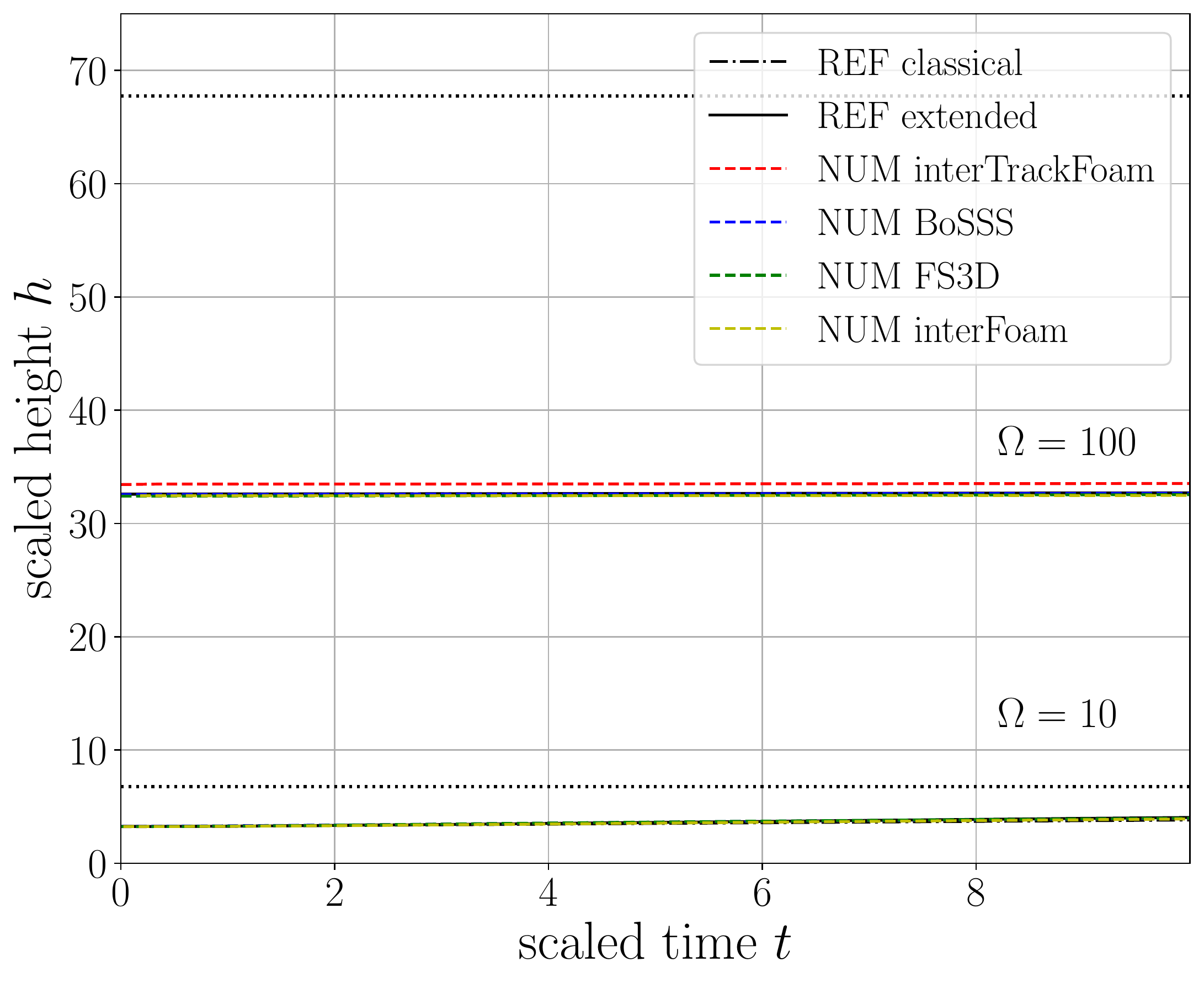}}
	\caption{Dynamical behaviour for the $\Omega$-study with scaling III. left: $\Omega = [0.1, 0.5, 1]$, right: $\Omega = [10, 100]$.} 
	\label{fig:omega_study_compar_Scale3}
\end{figure}

For completeness, Figure \ref{fig:omega_study_compar_Scale3} shows the results for $\Omega=10$ and $\Omega=100$. For scaling III, these curves are only horizontal lines, as they have been scaled such that the rise has ``barely started''. All solutions of the continuum mechanical problem basically agree with the reference solution in this initial state.

The overall results show the same quantitative behavior as the unscaled data. The agreement between the four numerical approaches is good. Furthermore, the numerical approaches are typically located ``between'' the two reference solutions. 

\newpage
\bibliography{references}              
\bibliographystyle{acm}
\end{document}